\documentclass[11pt]{amsart} 
\usepackage{amsthm,amsmath,amssymb,amscd} 
\usepackage{amssymb} 
\usepackage{graphics} 

\newtheorem{teor}{Theorem}[section] 
 
\newtheorem{defin}{Definition}[section] 
\newtheorem{remar}{Remark}[section] 
\newtheorem{prop}{Proposition}[section] 
\newtheorem{corol}{Corollary}[section] 
\newtheorem{lemma}{Lemma}[section]

\newcommand{\fdim}{\hfill {$\blacksquare$}}

\newcommand{\K}{K\"{a}hler}

\newcommand{\C}{\mathbb{C}} 
 
\newcommand{\e}{\varepsilon} 
\newcommand{\del}{\partial}

\begin{document} 

\title{Blowing up K\"{a}hler manifolds with constant 
scalar curvature II} 

\author[Claudio Arezzo] {Claudio Arezzo} 
\address{ claudio.arezzo@unipr.it \\ Dipartimento di Matematica\\ 
Universita' di Parma\\ 
Ê ÊVia M. D'Azeglio, 85\\43100 Parma\\Italy} 

\author[Frank Pacard] {Frank Pacard} 
\address{pacard@univ-paris12.fr\\ University Paris 12 and 
Institut Universitaire de France, France} 

\begin{abstract} 
In this paper we prove the existence of \K\ metrics of constant 
scalar curvature on the blow up at finitely many points of a compact 
manifold which already carries a \K\ constant scalar curvature 
metric. Necessary conditions of the number and locations of the blow 
up points are given. 
\end{abstract} 

\maketitle 

\vspace{-,15in} 


\section{Introduction} 

Assume that $(M, \omega_M)$ is a constant scalar curvature compact 
\K\ manifold. Given $m$ distinct points $p_1, \ldots, p_m \in M$, we 
define $\mbox{Bl}_{p_1, \ldots, p_m} M$ to be the blow up of $M$ at 
the points $p_1, \ldots, p_m$. The question we would like to address 
in the paper is whether $\mbox{Bl}_{p_1, \ldots, p_m} M$ can be 
endowed with a constant scalar curvature \K \ form. In 
\cite{Are-Pac}, we have already given a positive answer to this 
question in the case where the manifold has no nontrivial vanishing 
holomorphic vector field (i.e. has no non zero holomorphic vector 
field which vanish somewhere). This condition is for example 
fulfilled when the group of automorphisms of $M$ is discrete. Under 
this condition, we have obtained the following~: 
\begin{teor} \cite{Are-Pac} 
Assume that $(M, \omega_M)$ is a constant scalar curvature compact 
\K\ manifold without nontrivial vanishing holomorphic vector field. 
Then, $\tilde M$, the blow up at finitely many points has a constant 
scalar curvature \K\ form $\omega_{\tilde M}$. In addition, if the 
scalar curvature of $\omega_M$ is not zero then the scalar 
curvatures of $\omega_{\tilde M}$ and of $\omega_M$ have the same 
signs. 
\end{teor} 

In the case of $0$-scalar curvature metrics, we also have the~: 
\begin{teor} 
\label{sfk} \cite{Are-Pac} Assume that $(M, \omega_M)$ is a 
$0$-scalar curvature compact \K\ manifold without nontrivial 
vanishing holomorphic vector field. Then, the blow up of $M$ at 
finitely many points has a $0$-scalar curvature \K\ form provided 
the first Chern class of $M$ is non zero. 
\end{teor} 

This last result compliments, in any dimension, previous 
constructions which have been obtained in complex dimension $n=2$ 
and for $0$-scalar curvature metrics, by Kim-LeBrun- Pontecorvo 
\cite{klbp}, LeBrun-Singer \cite{lbs} and Rollin-Singer \cite{rs}. 
Indeed, using twistor theory, Kim-LeBrun- Pontecorvo and 
LeBrun-Singer have been able to construct such metrics by 
desingularizing some quotients of minimal ruled surfaces and, more 
recently, Rollin-Singer \cite{rs2} have shown that one can 
desingularize compact manifolds of zero scalar curvature with cyclic 
orbifold groups, keeping the scalar curvature zero by solving on the 
desingularization the hermitian anti-selfdual equation which implies 
the \K\ scalar flat equation. 

\medskip 

Theorem 1.1 and Theorem 1.2 are obtained using a connected sum of 
the \K\ form $\omega_M$ at each $p_\ell$ with a $0$-scalar curvature 
\K\ form $\eta$ which is defined on $N := \mbox{Bl}_0 \, {\mathbb 
C}^n$, the blow up of ${\mathbb C}^n$ at the origin. When $n=2$ this 
model $(N, \eta)$ has been known for long time and it is usually 
referred to in the literature as Burns metric \cite{lb}. In higher 
dimensions a similar model has been produced by Simanca \cite{si}. 

\medskip 

In the present paper, we are interested in the case where $M$ has 
nontrivial vanishing holomorphic vector fields $X_1, Ê\ldots, X_d$ , 
with $d \geq 1$ (this condition in particular implies that $M$ has a 
nontrivial continuous automorphisms group) since this is precisely 
the case which is not covered by our previous result. Our main 
result states that the blow up of $M$ at sufficiently many points 
can be endowed with a constant scalar curvature \K\ metric provided 
these points are carefully chosen. 
\begin{teor} 
\label{main1} Assume that $(M, \omega_M)$ is a constant scalar 
curvature compact \K\ manifold. Then, there exists $m_0 \geq Êd + 1$ 
and, for all $m \geq m_0$ there exist a nonempty open subset $U 
\subset M^m$ such that, for all $(p_1, \ldots, p_m) \in U$ there 
exist constant scalar curvature \K\ forms on the blow up of 
$\mbox{Bl}_{p_1, \ldots, p_m} \,M$. 
\end{teor} 

In fact, we can be more precise and show that there exists $\tilde 
\omega_\e$, for $\e \in (0, \e_0)$, a one parameter family of \K\ 
forms on Ê$\mbox{Bl}_{p_1, \ldots, p_m} \, M$ which converges in 
${\mathcal C}^\infty$ norm as $\e$ tends to $0$, to the \K\ form 
$\omega_M$ away from the points $p_\ell$. Moreover, the sequence of 
\K\ forms $\e^{-2} \, \tilde \omega_\e$ converges in ${\mathcal 
C}^\infty$ norm to the \K\ form $\eta$ on $N$ the blow up at the 
origin of ${\mathbb C}^n$. Finally, if $\omega_M$ had positive 
(resp. negative) scalar curvature then the \K\ forms $\tilde 
\omega_\e$ have positive (resp. negative) scalar curvature. 

\medskip 

Let us emphasize that, in contrast with the case where the manifold 
$M$ does not carry any nontrivial vanishing holomorphic vector 
field, this time the number and position of the blow up points is 
not arbitrary. 

\begin{remar} 
We will see that, in some explicit cases, one can make use of the 
symmetries of the manifold $M$, and construct constant scalar 
curvature \K\ forms on the blow up of $M$ at $(p_1, \ldots, p_m)$ 
with $m \leq d$. However this construction will not hold anymore for 
the choice of the points in some open set of $M^m$. 
\end{remar} 

Let us mention that, for the $2$-dimensional complex manifolds, \K\ 
metrics with constant scalar curvature have been obtained by 
Rollin-Singer in \cite{rs}, \cite{rs2} on blow ups of ${\mathbb 
P}^{1} \times {\mathbb P}^{1}$, ${\mathbb P}^{2}$ or ${\mathbb 
T}^{1} \times {\mathbb P}^{1}$ using a different approach based on 
both algebraic tools and a connect sum result, while for zero scalar 
curvature metrics on blow ups of ${\mathbb P}^{1} \times \Sigma_g$, 
$g\geq2$, the problem was solved by LeBrun-Singer \cite{lbs} using 
twistor theory. We will recall these results in \S $8$. 

\medskip 

In the statement of Theorem 1.3, the value of $m_0$ and the choice 
of the blow up points are very delicate issues. To describe this, we 
still assume that $(M, \omega_M)$ is a constant scalar curvature 
compact \K\ manifold and we define (see \S $3$) the elliptic fourth 
order operator 
\[ 
{\mathbb L}_M : = - \frac{1}{2}\, \Delta^{2}_M - \mbox{Ric}_M \cdot 
\nabla_M^{2}, 
\] 
where all operators are computed with respect to the metric induced 
by $\omega_M$. We assume that the kernel of ${\mathbb L}_M$ is not 
trivial and we denote by $\xi_0 \equiv 1, \xi_1, \ldots, \xi_d$ Êthe 
linearly independent functions which span this kernel. We can 
assume, without loss of generality that the functions $\xi_1, 
\ldots, \xi_d$ have mean $0$. Observe that $d$ either denotes the 
dimension of the space of vanishing holomorphic vector fields or the 
dimension of the nontrivial kernel of ${\mathbb L}_M$ since, by a 
result of LeBrun-Simanca \cite{ls} these two dimensions are the 
same. This being understood, we define the matrix 
\begin{equation} 
{\mathfrak M} (p_1, \ldots, p_m) : = \left( 
\begin{array}{cccclllll} 
\xi_1 (p_1) & \ldots & \xi_1 (p_m) \\[3mm] 
\vdots & Ê& \vdots \\[3mm] 
\xi_d(p_1) & \ldots & \xi_d (p_m) 
\end{array} \right) 
\label{eq:d-8} 
\end{equation} 
as well as the integer valued functions 
\[ 
{\mathfrak C}_1 : = \mbox{Rank} \, {\mathfrak M} \qquad \mbox{and} 
\qquad {\mathfrak C}_2 := \mbox{dim} \, (K^m_+ \cap \mbox{Ker} \, 
{\mathfrak M} ) \label{eq:d-9} 
\] 
where $K^m_+$ is the cone of vectors with positive entries in 
${\mathbb R}^m$. 

\medskip 

Our main result is a consequence of the following sequence of 
results. The first one gives a sufficient condition on the functions 
${\mathfrak C}_1 $ and ${\mathfrak C}_2 $ for Theorem 1.3 to hold~: 
\begin{prop} 
Assume that $(M, \omega_M)$ is a constant scalar curvature compact 
\K\ manifold. Let us assume that $m \geq 1$ and $p_1, \ldots, p_m$ 
are chosen so that ${\mathfrak C}_1 = d$ and that ${\mathfrak C}_2 
\neq 0 $. Then there exists $\e_0 > 0$ and for all $\e \in (0, 
\e_0)$, a constant scalar curvature \K\ form $\tilde \omega_\e$ on 
$\mbox{Bl}_{p_1, \ldots, p_m} \, M$. 
\end{prop} 
Furthermore, the \K\ form $\tilde \omega_\e$ satisfy the properties 
described after the statement of Theorem 1.3. 

\medskip 

As will be explained in \S $6$, the first condition ${\mathfrak C}_1 
= d$ is easily seen to be generic (and open) in the sense that~: 
\begin{lemma} 
With the above notations, assume that $m \geq d$ then the set of the 
points $(p_1, \ldots, p_m) \in M^m$ such that ${\mathfrak C}_1 = d$ 
is open and dense in $M^m$. 
\end{lemma} 
A result of LeBrun-Simanca \cite{ls} (see also \cite{Fu} and 
\cite{lbs}) implies that nontrivial vanishing holomorphic vector 
fields are in one to one correspondence with nonconstant elements of 
the kernel of ${\mathbb L}_M$. Furthermore, it is well known that, 
for a choice of blow up points $p_1, \ldots, p_d$ in some open and 
dense subset of $M^d$, the group of automorphisms of 
$\mbox{Bl}_{p_1, \ldots, p_d} \, M$ is trivial (observe that, by 
virtue of the above mentioned result of LeBrun and Simanca $d$ is 
also equal to the dimension of the identity component of the 
automorphisms group of $M$). In view of all these results, one is 
tempted to conjecture that our condition ${\mathfrak C}_1 = d$ is 
equivalent to the fact that the group of automorphisms of 
$\mbox{Bl}_{p_1, \ldots, p_{\tilde d}} \, M$ is trivial. However, 
these two conditions turn out to be of a different nature. Indeed, 
for example when $d=1$ the kernel of ${\mathbb L}_M$ is spanned by 
$\xi_0 \equiv 1$ and $\xi_1$, and the set of nontrivial vanishing 
holomorphic vector fields is spanned by $X_1$. Then it is enough to 
choose $p$ not in the zero set of $\xi_1$ to get ${\mathfrak C}_1 
=1$, while the group of automorphisms of $\mbox{Bl}_{p} \, M$ is 
trivial if and only if $p$ is chosen away from the $0$ set of $X_1$, 
but, according to LeBrun and Simanca, the $0$ set of $X_1$ 
corresponds to the set of critical points of $\xi_1$. 

\medskip 

The second sufficient condition on the existence of points for which 
${\mathfrak C}_2 \neq 0$ is more subtle. We prove, in \S $6$, that 
this condition is always fulfilled for some careful choice of the 
points, provided their number $m$ is chosen larger than some value 
$m_0 \geq d+1$. 
\begin{lemma} 
With the above notations, assume that $m \geq d+1$, then the set of 
points $(p_1, \ldots, p_m) \in M^m$ for which ${\mathfrak C}_1 Ê=d$ 
and ${\mathfrak C}_2 \neq 0$ is open. Moreover, there exists $m_0 
\geq d+1$ such that, for all $m \geq m_0$ the set of points $(p_1, 
\ldots, p_m) \in M^m$ for which ${\mathfrak C}_1 Ê=d$ and 
${\mathfrak C}_2 \neq 0$ is nonempty. 
\end{lemma} 
The proof of the Lemma 1.2 is due to E. Sandier, and we are very 
grateful to him for allowing us to present it here. By opposition to 
the first condition, it is easy to convince oneself that this 
condition does not hold for generic choice of the points. 

\medskip 

\begin{remar} 
We believe that $m_0$, the minimal number of points for which the 
second condition is fulfilled, should be equal to $d+1$. We have 
only been able to prove this fact for some special manifolds where 
all elements of the kernel of ${\mathbb L}_M$ are explicitly known. 
\end{remar} 

It is interesting to observe that our sufficient conditions depends 
only on the complex structure of the base manifold and not on 
riemannian properties of its \K\ metric. This is shared with all 
known necessary conditions, from Matsushima-Lichnerovicz's one, 
Futaki's character and the known stability conditions implied by the 
existence of canonical metrics thanks to the work of Donaldson, 
Mabuchi and Paul-Tian. 

\medskip 

We now give a number of explicit examples for which our result can 
be applied. As already mentioned in Remark 1.1, by definition $m_0$ 
is larger than the dimension of vanishing holomorphic vector fields 
on $M$ and, using some symmetries of the manifolds, one can still 
construct constant scalar curvature \K\ forms on the blow up of $M$ 
at fewer points. To illustrate this fact, we have, when $M={\mathbb 
P}^n$, we have~: 
\begin{lemma} 
When $M ={\mathbb P}^n$, then $d = n^{2} + 2n$ and $m_0 \leq Ê2 
n(n+1)$. 
\end{lemma} 
Applying Theorem 1.3, we get the existence of constant scalar 
curvature \K\ forms on the blow up of ${\mathbb P}^n$ at $m$ points 
for $m \geq m_0$ which belong to some nonempty open set of $M^m$. 
One can also work equivariantly using the symmetries of ${\mathbb 
P}^n$ to get the~: 
\begin{corol} 
For any $m \geq n+1$ there exist points $p_1,\dots, p_m$ such that 
the blow up of ${\mathbb P}^n$ at $p_1,\dots, p_m$ carries a 
constant scalar curvature \K\ metric. 
\end{corol} 

The above corollary is optimal in the number of points because, for 
$m\leq n$, $\mbox{Bl}_{p_1, \ldots, p_m} \, {\mathbb P}^n$ does not 
satisfy the Matsushima-Lichnerowicz obstruction. Observe also that 
$\mbox{Bl}_{p_1, \ldots, p_{n+1}} \, {\mathbb P}^n$ still has 
vanishing holomorphic vector fields. 

\medskip 

On the other hand it is well known that on ${\mathbb P}^n$, $n+2$ 
points forming a projective frame are enough to kill all holomorphic 
vector fields when blowing up, and indeed we can prove that this 
condition also guarantees the existence of \K\ constant scalar 
curvature metric. 
\begin{corol} 
Given $p_1, \dots , p_m$ points in ${\mathbb P}^n$ such that $p_1, 
\dots , p_{n+2}$ form a projective frame, $\tilde{M} = 
\mbox{Bl}_{p_1, \ldots, p_{m}}$ has a \K\ constant scalar curvature 
metric and no holomorphic vector fields. Moreover $p_1, \dots , 
p_{n+2}$ can vary in a dense open subset of $({\mathbb P}^n)^{n+2}$ 
and $p_{n+3}, \dots , p_m$ are arbitrary. 
\end{corol} 

When $M = {\mathbb P}^n \times M_0$ where $M_0$ is any manifold with 
no nontrivial vanishing holomorphic vector field, all the above 
result carry to $M$. In particular~: 
\begin{corol} 
Assume that $(M_0, \omega_0)$ is a \K\ manifold with constant scalar 
curvature and without nontrivial vanishing holomorphic vector fields 
and define $M ={\mathbb P}^n \times M_0$, then $d = n^{2} + 2n$ and 
$m_0 \leq 2n(n+1)$. Moreover, the blow up of ${\mathbb P}^n \times 
M_0 $ at $m \geq n+1$ points has a constant scalar curvature \K\ 
metric. 
\end{corol} 
As already mentioned, when $n=1$ and $M_0 = Ê\Sigma_g$ is a Riemann 
surface of genus $g$ greater than $1$, the existence of zero scalar 
curvature \K\ metrics on blow ups of ${\mathbb P}^1 \times \Sigma_g$ 
is due to LeBrun and Singer \cite{ls}. 

\medskip 

Finally, we have also considered the case where $M = {\mathbb P}^{n} 
\times {\mathbb P}^{m}$, for which we have obtained the~: 
\begin{corol} 
Suppose $n\leq m$. The blow up of ${\mathbb P}^{n} \times {\mathbb 
P}^{m}$ at $k \geq m+2$ special points has a constant scalar 
curvature \K\ metric. Moreover, if the projections of these points 
to the two factors contain projective frames, these points can move 
in a dense open subset of $({\mathbb P}^{n} \times {\mathbb 
P}^{m})^{m+2}$ and the other ones are arbitrary. 
\end{corol} 

The above results for $n=2$ translate directly in the following~: 
\begin{corol} 
For any $m \geq 2$ there exist points $p_1, \ldots , Êp_m \in 
{\mathbb P}^1 \times {\mathbb P}^1$ such that $\mbox{Bl}_{p_1, 
\ldots, p_m} \, {\mathbb P}^1 \times {\mathbb P}^1$ has a constant 
scalar curvature \K\ metric. 
\end{corol} 
Rollin-Singer in \cite{rs2} have already found such constant scalar 
curvature \K\ metrics on the blow up of ${\mathbb P}^1 \times 
{\mathbb P}^1$ at $m \geq 6$ points. It is worth mentioning that, 
for $m=2$, $\mbox{Bl}_{p_1, p_2} \, {\mathbb P}^1 \times {\mathbb 
P}^1$ still has nontrivial vanishing holomorphic vector fields. When 
$m \leq 7$ the construction is performed using the symmetries of the 
manifold, working equivariantly. 

\medskip 

\section{Weighted spaces} 

In this section, we describe weighted spaces on the noncompact 
manifold $M$ with $m$ points removed, as well as weighted spaces on 
the noncompact manifolds $N$. We also introduce notations which will 
be used in the next sections. 

\medskip 

For all $r > 0$, we agree that 
\begin{equation} 
\begin{array}{rlll} 
B_r & : = & \{ Êz \in {\mathbb C}^n \quad : \quad Ê|z| < r \}, \\[3mm] 
B^*_r & : = & \{ Êz \in {\mathbb C}^n Ê\quad : Ê\quad 0 < 
|z| \leq r \} ,\\[3mm] 
C_r & : = & \{ z \in {\mathbb C}^n Ê\quad : \quad |z | \geq r \} 
\end{array} 
\label{eq:f-5} 
\end{equation} 

Assume that $k \in {\mathbb N}$ and $\alpha \in (0,1)$ are fixed. 
Given $\bar r >0$ and a function $v \in {\mathcal C}^{k, 
\alpha}_{loc} (B^*_{\bar r})$, we define 
\[ 
\| v \|_{{{\mathcal C}}^{k, \alpha}_{\delta} (B^*_{\bar r})} : = 
\sup_{0 < r \leq \bar r} Êr^{-\delta} \, \| v (r \, \cdot) 
\|_{{{\mathcal C}}^{k, \alpha} (\bar B_1 -B_{1/2})}. 
\] 
and, for any function $v \in {\mathcal C}^{k, \alpha}_{loc} (C_{\bar 
r})$, we define 
\[ 
\| v \|_{{{\mathcal C}}^{k, \alpha}_{\delta} (C_{\bar r})} : = 
\sup_{r \geq \bar r} Êr^{-\delta} \, \| v (r \, \cdot) 
\|_{{{\mathcal C}}^{k, \alpha} ( \bar B_2 -B_1)} . 
\] 

Assume that we are given $(M, \omega_M)$, a K\"ahler manifold with 
K\"ahler form $\omega_M$ and that we are also given $m$ distinct 
points $p_1, \ldots, p_m \in M$. By definition, near $p_\ell$, the 
manifold $M$ is biholomorphic to a neighborhood of $0$ in ${\mathbb 
C}^n$. In particular, we can choose complex coordinates $z : = (z_1, 
\ldots, z_n)$ in a neighborhood of $0$ in ${\mathbb C}^n$, to 
parameterize a neighborhood of $p_\ell$ in $M$. In order to 
distinguish between the different neighborhoods and coordinate 
systems, we agree that, for all $r$ small enough, say $r \in (0, 
r_0)$, $B_{\ell , r}$ (resp. $B_{\ell , r}^*$) denotes the ball 
(resp. the punctured ball) of radius $r$ in the above defined 
coordinates $z$ parameterizing a fixed neighborhood $B_{\ell, r_0}$ 
of $p_\ell$. Furthermore, it follows from \cite{Gri-Har}, that there 
is no loss of generality in assuming that the \K\ form $\omega_M$ 
can be expended as 
\begin{equation} 
\omega_M : = i Ê\, Ê\del \, \bar \del \, (\mbox{$\frac{1}{2}$} \, 
|z|^2 + \varphi_\ell ), \label{eq:normal} 
\end{equation} 
where the function $\varphi_{ \ell} Ê\in {\mathcal C}^{3, \alpha}_4 
(B^*_{\ell , r_0})$. This in particular implies that, in these 
coordinates, the standard metric on ${\mathbb C}^n$ and the metric 
induced by $\omega_M$ agree up to order $2$. 

\medskip 

For all $r \in (0, r_0)$, we define 
\begin{equation} 
M_r : = M - \cup_\ell \, B_{\ell , r} . \label{eq:f-6} 
\end{equation} 
The weighted space for functions defined on the noncompact manifold 
\begin{equation} 
M^* : = M - \{ p_\ell , \quad : \quad \ell=1, \ldots, m Ê\}. 
\label{eq:f-4} 
\end{equation} 
is then defined as the set of functions whose decay or blow up near 
any $p_\ell$ is controlled by a power of the distance to $p_\ell$. 
More precisely, we have the~: 
\begin{defin} 
Given $k \in {\mathbb N}$, $\alpha \in (0,1)$ and $\delta \in 
{\mathbb R}$, we define the weighted space ${\mathcal C}^{k, 
\alpha}_{\delta} (M^*)$ to be the space of functions $w \in 
{{\mathcal C}}^{k, \alpha}_{loc} (M^*)$ for which the following norm 
is finite 
\[ 
\| w \|_{{{\mathcal C}}^{k, \alpha}_{\delta} (M^*)} Ê: = \| w 
\|_{{{\mathcal C}}^{k, \alpha} (M_{r_0/2})} + \sup_\ell \, \| w 
|_{B_{\ell , r_0}^*} \|_{{{\mathcal C}}^{k, \alpha}_\delta (B_{\ell 
,r_0}^* )}. 
\] 
\label{de:f-4.1} 
\end{defin} 

We now turn to the description of weighted space on $(N, \eta)$, the 
blow up at the origin of ${\C}^n$ endowed with the Burns-Simanca 
metric. Away from the exceptional divisor, the \K \ form $\eta$ is 
given by 
\[ 
\eta Ê= i \, \del \, \bar \del f_n (|v|^2) 
\] 
where $v = (v_1, \ldots, v_n)$ are complex coordinates in ${\mathbb 
C}^n -\{ 0 \}$ and where the function $s \longrightarrow f_n(s)$ is 
a solution of the ordinary differential equation 
\[ 
s2 \, ( s\, \del_s f_n )^{n-1} \, \del_s2 f_n + (n-1) \, Ê s \, 
\del_s f_n Ê- Ê(n-2) =0 
\] 
which satisfies $f_n \sim \log s$ near $0$. It turns out that, when 
$n=2$, the function $P_2$ is explicitly given by 
\[ 
f_2 (s) = Ê\log Ês + \lambda \, s + c 
\] 
where $\lambda >0$ and $c \in {\mathbb R}$, while in dimension $n 
\geq 3$, even though there is no explicit formula for $f_n$ we have 
the following~: 
\begin{lemma} 
Assume that $n \geq 3$. Then the function $f_n$ can be expanded as 
\[ 
f_n (s) = \lambda \, s + c - \lambda^{2-n} \, s^{2-n} + {\mathcal O} 
(s^{1-n}) 
\] 
for $s >1$, where $\lambda >0$ and $c \in {\mathbb R}$. 
\end{lemma} 
{\bf Proof~:} This result follows from the analysis done in 
\cite{si} even though this is not explicitly stated. In fact, we 
define the function $\zeta$ by $ s \, \zeta : = s \, \del_s f_n - 
1$. Direct computation shows that $\zeta$ solves 
\[ 
(1 + s \, \zeta)^{n-1} \, s^{2} \, \del_s \zeta = Ê(1+ s\, 
\zeta)^{n-1} - 1 - (n-1) \, s \, \zeta 
\] 
If in addition we take $\zeta (0) = 1$, then $\del_s \zeta $ remains 
positive and one can check that $\zeta$ is well defined for all time 
and converges to some positive constant $\lambda$, as $s$ tends to 
$\infty$. This immediately implies that $s \, \del_s f_n = \lambda 
\, s + {\mathcal O} (1)$ at infinity. The expansion then follows 
easily. \fdim 

\medskip 

Changing variables $u : = Ê\sqrt{2 Ê\, \lambda} \, v$, we see from 
the previous Lemma that the \K\ form $\eta$ can be expanded near 
$\infty$ as 
\begin{equation} 
\eta = i \, \del \, \bar \del ( \mbox{$\frac{1}{2}$} \, |u|^2 + \log 
|u|^2 + \tilde \varphi ) , \label{eq:f-111} 
\end{equation} 
in dimension $n =2$ and as 
\begin{equation} 
\eta = i \, \del \, \bar \del ( \mbox{$\frac{1}{2}$} \, |u|^2 - 
2^{n-2} \, |u|^{4-2n} + \tilde \varphi) , \label{eq:f-1111} 
\end{equation} 
in dimension $n \geq 3$, where $\tilde \varphi \in {\mathcal C}^{3, 
\alpha}_{3-2n} (C_1)$ for all $n \geq 2$. The property which will be 
crucial for our construction is that, in these two expansions, the 
coefficient in from of $\log |u|^2$, in dimension $n=2$ and the 
coefficients in front of $|u|^{4-2n}$ are not zero. 

\medskip 

For all $R > 1$, we define 
\begin{equation} N_{R} Ê: = N - C_R. 
\label{eq:f-8} 
\end{equation} 

We are now in a position to define weighted spaces on the noncompact 
complete manifold $N$. This time, we are interested in functions 
which decay (or blow up) near infinity at a rate which is controlled 
by a power of the distance to a fixed point in $N$. More precisely, 
we have the~: 
\begin{defin} 
Given $k \in {\mathbb N}$, $\alpha \in (0,1)$ and $\delta \in 
{\mathbb R}$, we define the weighted space ${{\mathcal C}}^{k, 
\alpha}_{\delta} (N)$ to be the space of functions $w \in {{\mathcal 
C}}^{k, \alpha}_{loc} (N)$ for which the following norm is finite 
\[ 
\| w \|_{{{\mathcal C}}^{k, \alpha}_{\delta} (N)} Ê: = \| w 
\|_{{{\mathcal C}}^{k, \alpha} (N_{2})} + Ê\| w |_{C_{1}} 
\|_{{{\mathcal C}}^{k, \alpha}_\delta (C_{1})} . 
\] 
\label{de:4.2} 
\end{defin} 

\section{The geometry of the equation} 

Recall that $(M,\omega_M)$ is a compact \K\ manifold of complex 
dimension $n$. We will indicate by $g_M$ the Riemannian metric 
induced by $\omega_M$, $\mbox{Ric}_M$ its Ricci tensor, $\rho_M$ the 
corresponding Ricci form, and ${\bf s} (\omega_M) = s_M$ its scalar 
curvature. 

\medskip 

Following LeBrun-Simanca \cite{ls}, we want to understand the 
behavior of the scalar curvature of the metric under deformations of 
the form 
\[ 
\omega : = \omega_M Ê+ i \, \del \, \bar \del \, \varphi 
\] 
where $\varphi$ a function defined on $M$. The following result is 
proven, for example in \cite{ls}~: 
\begin{prop} 
The scalar curvature of $\omega $ can be expanded in terms of 
$\varphi$ as 
\[ 
{\bf s} Ê\, (\omega ) = {\bf s} \, Ê(\omega_M ) - \frac{1}{2} \, 
\Delta^{2}_M \, \varphi - \mbox{Ric}_M \cdot \nabla_M^{2} \, \varphi 
+ Q_M (\varphi), 
\] 
where $Q_M$ collects all the nonlinear terms and where all operators 
in the right hand side of this identity are computed with respect to 
the \K\ form $\omega_M$. \label{pr:f-2.1} 
\end{prop} 

The operator 
\begin{equation} 
{\mathbb L}_M Ê: = -\frac{1}{2} \, \Delta^{2}_M Ê- \mbox{Ric}_M 
\cdot \nabla_M^{2}, \label{eq:f-2} 
\end{equation} 
which acts on functions is the one which has been defined in the 
introduction. 

\medskip 

For a general \K\ metric it can be very difficult to analyze these 
operators. Nevertheless geometry comes to the rescue at a constant 
scalar curvature metric. Indeed, in this case we have 
\begin{equation} 
{\mathbb L}_M Ê = ({\bar{\partial}} 
\partial^{\#})^{\star} {\bar{\partial}}\partial^{\#}, 
\label{eq:f-3} 
\end{equation} 
where $\partial^{\#} \, \varphi $ denotes the $(1,0)$-part of the 
$g_M$-gradient of $\varphi$. Using this result, the key observation 
of LeBrun-Simanca is that to any element $\varphi$ of $\mbox{Ker} \, 
{\mathbb L}_M$ one can associate a holomorphic vector field, namely 
$\partial^{\#} \varphi $, which must vanish somewhere on $M$. This 
is the crucial relation between constant scalar curvature (or 
extremal) metrics and the space of holomorphic vector fields $H(M)$ 
(and, in turn, the automorphism group $\mbox{Aut} \, (M)$). 

\medskip 

We recall the following important result, which has been mentioned 
in the introduction and which helps to understand the kernel of the 
operator ${\mathbb L}_M$~: 
\begin{teor} 
\label{clsa-1} \cite{ls} Assume that $(M, \omega_M)$ is a compact 
\K\ manifold of complex dimension $n$, then the space of the 
elements of the kernel of ${\mathbb L}_M$ whose mean over $M$ is 
$0$, is in one to one correspondence with the space of holomorphic 
vector fields which vanish on $M$. \label{th:f-2.1} 
\end{teor} 

The previous considerations extend to $(N, \eta)$ and this implies 
the following important result of \cite{Kov-Sin} which states that 
there are no elements in the kernel of the operator ${\mathbb 
L}_{N}$ which decay at infinity. 
\begin{prop} 
\label{clsa-11} \cite{Kov-Sin} There are no nontrivial solution to 
${\mathbb L}_N \, \varphi = 0$, which belong to ${\mathcal C}^{4, 
\alpha}_{\delta} (N)$, for some $\delta <0$. 
\end{prop} 
{\bf Proof~:} Assume that ${\mathbb L}_N \, \varphi =0$ and that 
$\varphi \in{\mathcal C}^{4, \alpha}_{\delta} (N)$, for some $\delta 
<0$. Then $\del^\# \varphi$ is an holomorphic vector field which 
tends to $0$ at infinity. Using Hartogs theorem, the restriction of 
$\del^\# \varphi $ to $C_1$ can be extended to a holomorphic vector 
field on ${\mathbb C }^n$. Since this vector field decays at 
infinity, it has to be identically equal to $0$. This implies that 
$\del^\# \varphi$ is identically equal to $0$ on $C_1$. However 
$\varphi$ being a real valued function, this implies that $\del 
\varphi = \bar \del \varphi = 0$ in $C_1$. Hence the function 
$\varphi$ is constant in $C_1$ and decays at infinity. This implies 
that $\varphi$ is identically equal to $0$ in $C_1$ and satisfies 
${\mathbb L}_N \, \varphi =0$ in $N$. Now, use the unique 
continuation theorem to conclude that $\varphi$ is identically equal 
to $0$ in $N$. \fdim 

\section{Mapping properties} 

\subsection{Analysis of the operator defined on $M^*$} 

The results we want to obtain are based on the fact that, near each 
$p_\ell$ and in the above defined coordinates, the metric on $M$ is 
asymptotic to the Euclidean metric. This implies that, in each 
$B_{\ell , r_0 }$ the operator ${\mathbb L}_M$ is close to the 
operator 
\begin{equation} 
{\mathbb L}_0 : = - \frac{1}{2} \, \Delta_0^{2}, \label{eq:f-12} 
\end{equation} 
where $\Delta_0$ denotes the Laplacian in ${\mathbb C}^n$ when 
endowed with the standard \K\ form. In particular these two 
operators have the same indicial roots at $p_\ell$. This set of 
indicial roots is given by ${\mathbb Z} -\{ 5 - 2n, \ldots, -1\}$ 
when $n\geq 3$ and is given by ${\mathbb Z}$ when $n=2$. The mapping 
properties of ${\mathbb L}_M$ when defined between weighed spaces 
follows from the general theory explained in \cite{Mel} and this 
will allow us to construct some right inverse for ${\mathbb L}_M$. 

\medskip 

We define in $B^*_{\ell , r_0}$ the function $G_\ell$ by 
\[ 
G_\ell (z) = Ê- \log \, |z|^2 \qquad \mbox{ when $n=2$ \qquad and} 
\qquad G_\ell (z) = |z|^{4-2n} \qquad \mbox{ when $n\geq 3$}. 
\] 
Observe that, unless the metric $g_M$ is the Euclidean metric, these 
functions are not solutions of the homogeneous equation associated 
to ${\mathbb L}_M$, however they can be perturbed into $\tilde 
G_\ell$ solutions of the homogeneous problem ${\mathbb L}_M \, 
\tilde G_\ell =0$. Indeed, reducing $r_0$ if this is necessary, we 
have the following~: 
\begin{lemma} 
There exist functions $\tilde G_\ell$ which are solutions of 
${\mathbb L}_M \, \tilde G_\ell =0$ in $B_\ell^* \, (r_0)$ and which 
are asymptotic to $G_\ell$ in the sense that $\tilde G_\ell- G_\ell 
\in {\mathcal C}^{4, \alpha}_{6-2n} (B^*_{\ell ,r_0})$ when $n \geq 
4$ and $\tilde G_\ell- G_\ell \in {\mathcal C}^{4, \alpha}_{\delta} 
(B^*_{\ell ,r_0})$ for any $\delta <6-2n$ when $n=2,3$. 
\end{lemma} 
{\bf Proof~:} Observe that ${\mathbb L}_M \, G_\ell = ({\mathbb L}_M 
- {\mathbb L}_0) \, G_\ell$ and, thanks to the expansion given in 
(\ref{eq:normal}), we conclude that ${\mathbb L}_M \, G_\ell \in 
{\mathcal C}^{0, \alpha}_{2-2n} (B_{\ell Ê, r_0}^*)$. The result 
then follows from the analysis of \S 6 in \cite{Are-Pac}. Observe 
that, as long as $6-2n$ is not an indicial root of ${\mathbb L}_M$ 
we conclude that $\tilde G_\ell - G_\ell \in {\mathcal C}^{4, 
\alpha}_{6-2n} (B_{\ell , r_0}^*)$ while, if $6-2n$ is an indicial 
root of ${\mathbb L}_M$, we only conclude that $\tilde G_\ell- 
G_\ell \in {\mathcal C}^{4, \alpha}_{\delta} (B^*_{\ell ,r_0})$ for 
any $\delta <6-2n$. \fdim 

\medskip 

Observe that there is no uniqueness of $\tilde G_\ell$. With the 
functions $\tilde G_\ell$ at hand, we define the following 
deficiency spaces 
\[ 
{\mathcal D}_0 : = \mbox{Span} \{ \chi_1 Ê, Ê\ldots, \chi_m Ê\} , 
\qquad \mbox{and} \qquad {\mathcal D}_1 : = \mbox{Span} \{ \chi_1 \, 
\tilde G_1 , \ldots, \chi_m \, \tilde G_m\} , 
\] 
where $\chi_\ell$ is a cutoff function which is identically equal to 
$1$ in $B_{\ell ,r_0/2}$ and identically equal to $0$ in $M-B_{\ell 
, r_0 }$. 

\medskip 

When $n \geq 3$, we fix $\delta \in (4-2n, 0)$ Êand define the 
operator 
\[ 
\begin{array}{rclcllll} 
L_\delta : & Ê( {{\mathcal C}}^{4, \alpha}_{\delta} (M^*) \oplus 
{\mathcal D}_1 ) \times {\mathbb R} & \longrightarrow & Ê{{\mathcal 
C}}^{0,\alpha}_{\delta - 4} (M^*) \\[3mm] 
& (\varphi , \beta) Ê& \longmapsto & {\mathbb L}_M \, \varphi + 
\beta , 
\end{array} 
\] 
Whereas, when $n=2$, we fix $\delta \in (0,1)$ and define the 
operator 
\[ 
\begin{array}{rclcllll} 
L_\delta : & ({{\mathcal C}}^{4, \alpha}_{\delta} (M^*) \oplus 
{\mathcal D}_0 \oplus {\mathcal D}_1 ) \times {\mathbb R} & 
\longrightarrow & {{\mathcal C}}^{0,\alpha}_{\delta - 4} (M^*) \\[3mm] 
& (\varphi , \beta) & \longmapsto & {\mathbb L}_M \, \varphi + \beta 
, 
\end{array} 
\] 
To keep notations short, it will be convenient to set ${\mathcal D} 
: = {\mathcal D}_1$ when $n \geq 3$ and ${\mathcal D} : = {\mathcal 
D}_0 \oplus {\mathcal D}_1$ when $n=2$. The main result of this 
section reads~: 
\begin{prop} 
Assume that the points $p_1, \ldots, p_m \in M$ are chosen so that 
${\mathfrak C}_1 = d$. Then, the operator $L_\delta$ defined above 
is surjective (and has a $(m+1)$-dimensional kernel). 
\label{pr:f-5.1} 
\end{prop} 
{\bf Proof~:} The proof of this result follows from the general 
theorem described in \cite{Mel}, however, we choose here to describe 
an almost self contained proof. Recall that the kernel of ${\mathbb 
L}_M$ is spanned by the functions $\xi_0 \equiv 1, \xi_1, \ldots, 
\xi_d$ and, by assumption $\xi_j$, for $j=1, \ldots, d$ have mean 
$0$. We use the fact that ${\mathbb L}_M$ is self adjoint and hence, 
for $\psi \in L1 (M)$, the problem 
\[ 
{\mathbb L}_M \, \varphi = \psi 
\] 
is solvable if and only if $\psi$ satisfies 
\[ 
\int_M \psi \, \xi_j =0 
\] 
for $j=0, \ldots, d$. 

\medskip 

Observe that ${\mathcal C}^{0, \alpha}_{\delta-4} (M^*) \subset L1 
(M)$ when $\delta > 4-2n$. Now, given $\psi \in L1 (M)$, we choose 
\[ 
a_0 = \oint_M \psi 
\] 
and, since ${\mathfrak C}_1 = d$, we also choose $a_1, \ldots, a_d 
\in {\mathbb R}$ solution of the system 
\[ 
\forall j=1, \ldots, d , \qquad \qquad \int_M \psi \, \xi_j = 
\sum_{\ell = 1}^m a_\ell \, \xi_j (p_\ell) 
\] 
Then, the problem 
\[ 
{\mathbb L}_M \, \varphi = \psi - a_0 - \sum_{\ell=1}^m a_\ell \, 
\delta_{p_\ell} 
\] 
is solvable in $L1(M)$ for all $p \in [1, \frac{n}{n-1})$ and 
uniqueness follows if we impose that $\varphi$ is orthogonal to 
$\xi_0, \ldots, \xi_d$. To complete the proof, we invoke regularity 
theory which implies that $\varphi \in {{\mathcal C}}^{4, 
\alpha}_{\delta} (M^*) \oplus {\mathcal D}_1$ when $n \geq 3$ and 
${{\mathcal C}}^{4, \alpha}_{\delta} (M^*) \oplus {\mathcal D}_0 
\oplus {\mathcal D}_1$ when $n=2$. The estimate of the dimension of 
the kernel is left to the reader since it will not be used in the 
paper. \fdim 

\medskip 

Observe that, when we will solve ${\mathbb L} \, \varphi + \beta = 
\psi$, the constant $\alpha$ is determined by 
\[ 
\beta = \oint_M \psi. 
\] 

\subsection{Analysis of the operator defined on $N$} 

As above, we use the fact that the metric on $N$ is asymptotic to 
Euclidean metrics, as explained in (\ref{eq:f-111}) and 
(\ref{eq:f-1111}). This implies that, in $C_1$, the operator 
${\mathbb L}_{N}$ is close to the operator ${\mathbb L}_0$ and they 
have the same set of indicial roots which is given by ${\mathbb Z} 
-\{ 5 - 2n, \ldots, -1\}$ when $n\geq 3$ and is given by ${\mathbb 
Z}$ when $n=2$. Given $\delta \in {\mathbb R}$, we define the 
operator 
\[ 
\begin{array}{rclcllll} 
\tilde L_\delta : & {{\mathcal C}}^{4, \alpha}_{\delta} (N) & 
\longrightarrow & Ê{{\mathcal C}}^{0, 
\alpha}_{\delta - 4} (N) Ê\\[3mm] 
& \varphi & \longmapsto & {\mathbb L}_{N} \, \varphi 
\end{array} 
\] 
and recall the following result from \cite{Are-Pac} 
\begin{prop} 
Assume that $\delta \in (0,1)$. Then the operator $\tilde L_\delta$ 
defined above is surjective and has a one dimensional kernel spanned 
by the constant function. \label{pr:f-5.4} 
\end{prop} 
{\bf Proof~:} ÊThe result of Proposition~\ref{clsa-11} precisely 
states that the operator $\tilde L_{\delta'}$ is injective when 
$\delta' <0$. This implies that the operator $\tilde L_\delta$ is 
surjective when $\delta > 4-2n$. When $\delta \in (0,1)$, this also 
implies that the operator $\tilde L_\delta$ has a one dimensional 
kernel, spanned by the constant function. \fdim 

\subsection{Bi-harmonic extensions} 

We end up this section by the following simple result which follows 
at once from the application of the maximum principle and whose 
proof can be found in \cite{Are-Pac}. 
\begin{prop} 
Assume that $n \geq 2$. Given $h \in {\mathcal C}^{4, \alpha} (\del 
B_1 )$, $k \in {\mathcal C}^{2, \alpha} (\del B_1 )$ there exists a 
function $H^i_{h,k} \in {\mathcal C}^{4, \alpha} (B_1 )$ such that 
\[ 
{\mathbb L}_0 \, ÊH^i_{h,k} = 0 \qquad \mbox{in} \qquad B_1 
\] 
with 
\[ 
H^i_{h,k} = h \quad \mbox{and} \quad \Delta_0 H^i_{h,k} = k \qquad 
\mbox{on} \qquad \del B_1 Ê. 
\] 
Moreover, 
\[ 
\| H^i_{h,k} \|_{{\mathcal C}^{4, \alpha} (B_1 ) } \leq c \, 
(\|h\|_{{\mathcal C}^{4, \alpha}(\del B_1 )} + \|k \|_{{\mathcal 
C}^{2, \alpha}(\del B_1)}) 
\] 
\label{pr:f-5.5} 
\end{prop} 

We will also need the following result which differs slightly from 
what we have already used in \cite{Are-Pac}. The rational being that 
there exists a bi-harmonic extension of the data $(h,k)$ which is 
defined on the complement of the unit ball and decays at infinity. 
Moreover, this function is bounded by a constant times the distance 
to the origin to the power $4-2n$ (when $n \geq 3$). In the case 
where the function $k$ is assumed to have mean $0$, then the rate of 
decay can be improved and estimated as the distance to the origin to 
the power $3-2n$~: 
\begin{prop} Assume that $n \geq 2$. Given $h \in {\mathcal 
C}^{4, \alpha} (\del B_1 )$, $k \in {\mathcal C}^{2, \alpha} (\del 
B_1 )$ such that 
\[ 
\int_{\del B_1 } k =0 
\] 
there exists a function $H^o_{h,k} \in {\mathcal C}^{4 , 
\alpha}_{3-2n} (C_1)$ such that 
\[ 
{\mathbb L}_0 \, H^o_{h,k}= 0, \qquad \mbox{in} \qquad C_1 
\] 
with 
\[ 
H^o_{h} = h \quad \mbox{and} \quad \Delta_0 H^o_{h,k} =k \qquad 
\mbox{on} \qquad \del B_1 . 
\] 
Moreover, 
\[ 
\| H^o_{h,k} \|_{{\mathcal C}^{4, \alpha}_{3-2n} (C_1 ) } \leq c \, 
(\|h\|_{{\mathcal C}^{4, \alpha}(\del B_1 )} + \|k \|_{{\mathcal 
C}^{2, \alpha}(\del B_1)}) 
\] 
\label{pr:f-5.5bis} 
\end{prop} 

\subsection{Perturbation of $\omega_M$} 

We consider the \K\ form 
\begin{equation} 
\omega = \omega_M + i \, \del \, \bar \del \, \varphi . 
\label{eq:f-17} 
\end{equation} 
The scalar curvature of $\omega$ is given by 
\[ 
{\bf s} (\omega) = s (\omega_M) + Ê{\mathbb L}_M \, \varphi + Q_M 
(\varphi ) , 
\] 
where the operator ${\mathbb L}_M$ has been defined in 
(\ref{eq:f-2}) and where $Q_M$ collects all the nonlinear terms. The 
structure of $Q_M$ is quite complicated however, it follows from the 
explicit computation of the Ricci curvature that, near $p_\ell$, the 
nonlinear operator $Q_M$ can be decomposed as 
\begin{equation} 
Q_M (\varphi) = Ê\sum_j B_{\ell,j} (\nabla^{4} \varphi , \nabla^{2} 
\varphi) \, Q_{\ell,j} (\nabla^{2} \varphi) + \sum_{j'} \tilde 
B_{\ell, j'} (\nabla^{3} \varphi , \nabla^{3} \varphi) \, \tilde 
Q_{\ell, j'} (\nabla^{2} \varphi) ,\label{eq:f-14} 
\end{equation} 
where the $B_{\ell,j}$ and $\tilde B_{\ell ,j'}$ are bilinear forms 
with bounded coefficients and the sums are finite. In addition, 
identifying the Hessian of $\varphi$ with a symmetric matrix, one 
can check that the nonlinear operators $Q_{\ell, j}$ and $\tilde 
Q_{\ell, j'}$ are uniformly bounded and uniformly Lipschitz over the 
set of symmetric matrices whose entries are bounded by some fixed 
constant $c_0$. More precisely, they satisfy 
\begin{equation} 
| Q_{\ell,j} (A) | + | \tilde Q_{\ell,j'} (A) | Ê\leq c , 
\label{eq:f-15} 
\end{equation} 
and 
\begin{equation} 
| Q_{\ell,j} (A) - Q_{\ell,j} (A')| \leq c \, |A-A'| , 
\label{eq:f-16} 
\end{equation} 
for all symmetric matrices $A$ and $A'$ satisfying $|A| \leq c_0$ 
and $|A'|\leq c_0$. 

\medskip 

As usual, we denote by Ê$\xi_0, \ldots, \xi_d$, independent 
functions which span the kernel of ${\mathbb L}_M$. We agree that 
$\xi_0\equiv 1$ is the constant function and that all functions 
$\xi_1, \ldots, \xi_d$ are normalized to have mean $0$. Now, assume 
that we are given $a_0, a_1, \ldots, a_m > 0$ such that there exists 
a solution of 
\[ 
{\mathbb L}_M \, H_a = a_0 - c_n \, \sum_\ell a_\ell \, 
\delta_{p_\ell} 
\] 
where the constant $c_n$ is defined by 
\[ 
c_n := 4 (n-2)(n-1) \, |S^{2n-1}| \qquad \mbox{ when $n \geq 3$\quad 
and} \qquad c_2 : = 2 \, |S^{3}| 
\] 
Observe that such a function $H_a$ exists if and only if $a_0$ is 
given by 
\[ 
a_0 = c_n \sum_{\ell=1}^m a_\ell 
\] 
and the coefficients $a_1, \ldots, a_m $ are solutions of the system 
\[ 
\sum_{\ell=1}^m Êa_\ell \, \xi_j (p_\ell) = 0 
\] 
for $j=1, \ldots, d$. This amounts to ask that the vector $(a_1, 
\ldots, a_m)$ has positive entries and is in the kernel of the 
matrix ${\mathfrak M}(p_1, \ldots, p_m)$. This is precisely the 
origin of the condition ${\mathfrak C}_2 \neq 0$. It is not hard to 
check that~: 
\begin{lemma} 
Near each $p_\ell$, the function $H_a$ satisfies 
\[ 
H_a + a_\ell \, G_\ell + b_\ell \in Ê{\mathcal C}^{4, \alpha}_{1} 
(B_{\ell, r_0}^*) 
\] 
for some constant $b_\ell \in {\mathbb R}$. 
\end{lemma} 

We fix 
\begin{equation} 
r_\e: = \e^{\frac{2n-1}{2n+1}} \label{eq:re} 
\end{equation} 

Assume that we are given $h_\ell \in {\mathcal C}^{4, \alpha}(\del 
B_1)$ and $k_\ell \in {\mathcal C}^{2, \alpha}(\del B_1)$, $\ell=1, 
\ldots, m$, satisfying 
\begin{equation} 
\| h_\ell\|_{{\mathcal C}^{4, \alpha} (\del B_1)} + \| 
k_\ell\|_{{\mathcal C}^{2, \alpha} (\del B_1)} \leq \kappa \, 
r_\e^{4} ,\label{eq:f-18} 
\end{equation} 
where $\kappa >0$ will be fixed later on. Further assume that 
\begin{equation} \int_{\del B_1} k_\ell = 0 \label{eq:orth} 
\end{equation} 
We define 
\begin{equation} 
H_{h,k} : = Ê\sum_\ell \chi_\ell \, H^o_{h_\ell, k_\ell} ( \cdot / 
r_\e) , \label{eq:f-19} 
\end{equation} 
where $h$ and $k$ stand for $(h_1, \ldots, h_m)$ and $(k_1, \ldots, 
k_m)$ and where the cutoff functions $\chi_\ell$ are identically 
equal to $1$ in $B_{\ell,r_0/2}$ and identically equal to $0$ in $M 
- B_{\ell ,r_0}$. 

\medskip 

We fix $\delta \in (4-2n, 5-2n)$. The space ${\mathcal C}^{4, 
\alpha}_{\delta} (M_{r_\e}) \oplus {\mathcal D}$ (resp. $ {\mathcal 
C}^{0, \alpha}_{\delta-4} (M_{r_\e})$) is just defined as the space 
of restrictions of functions in ${\mathcal C}^{4, \alpha}_{\delta} 
(M^*) \oplus {\mathcal D}$ (resp. ${\mathcal C}^{0, 
\alpha}_{\delta-4} (M^*)$) to $M_{r_\e}$, endowed with the induced 
norm. 

\medskip 

We would like to solve the equation 
\[ 
{\bf s} \left( \omega_M + i \, \del \bar \del (\e^{2n-2} \, H_a + 
H_{h,k} + \varphi) \right) = {\bf s} \, (\omega_M) - \e^{2n-2} \, 
a_0 + \beta , 
\] 
where $\varphi \in {\mathcal C}^{4, \alpha}_{\delta} (M_{r_\e}) 
\oplus {\mathcal D}$ and $\beta \in {\mathbb R}$ have to be 
determined. This amounts to solve the equation 
\begin{equation} 
{\mathbb L}_M \, (H_{h,k} + \varphi) + Q_M (\e^{2n-2} \, H_a + 
H_{h,k} + \varphi ) = \beta , \label{eq:f-21} 
\end{equation} 
in $M_{r_\e}$. We consider an extension (linear) operator 
\[ 
{\mathcal E}_\e : {\mathcal C}^{0, \alpha}_{\delta-4} (M_{r_\e}) 
\longrightarrow {\mathcal C}^{0, \alpha}_{\delta-4} (M^*) , 
\] 
such that ${\mathcal E}_\e \, \psi = \psi$ in $M_{r_\e/2}$, 
${\mathcal E}_\e \, \psi$ is supported in $M_{r_\e/2}$ and $\| 
{\mathcal E}_\e\|\leq 2$. 

\medskip 

The equation we would like to solve can be rewritten as 
\begin{equation} 
{\mathbb L}_M \, \varphi Ê= - Ê{\mathcal E}_\e \left( \left( 
{\mathbb L}_M \, H_{h,k} + Q_M ( \e^{2n-2} \, H_a + H_{h,k} + 
\varphi Ê) \right)|_{M_{r_\e}} \right) + \beta . \label{eq:f-22} 
\end{equation} 
where $\beta$ is determined by 
\[ 
\beta = Ê\oint {\mathcal E}_\e \left( \left( {\mathbb L}_M \, 
H_{h,k} + Q_M ( \e^{2n-2} \, H_a + H_{h,k} + \varphi Ê) 
\right)|_{M_{r_\e}} \right) 
\] 

Granted these notations, we prove the~: 
\begin{lemma} 
We fix $\delta Ê\in (4-2n , 5-2n)$. There exists $c >0$, $\bar 
c_\kappa > 0$ and there exists $\e_0 >0$ such that, for all $\e \in 
(0, \e_0)$ 
\[ 
\| {\mathbb L}_M \, H_{h,k} \|_{{\mathcal C}^{0, \alpha}_{\delta 
-4}(M_{r_\e})} \leq \bar c_\kappa \, r_\e^{2n+1} , 
\] 
and 
\[ 
\| Q_M (\e^{2n-2} \, H_a + H_{h,k} ) \|_{{\mathcal C}^{0, 
\alpha}_{\delta -4}(M_{r_\e})} \leq c Ê\, \e^{4n-4} \, r_\e^{6 -4n 
-\delta } , 
\] 
where the norms are the restrictions of the weighted norms in 
Definition~\ref{de:f-4.1} to functions defined in $M_{r_\e}$. In 
addition, we have and 
\[ 
\int_{M_{r_\e}} \left| {\mathbb L}_M \, H_{h,k} + Q_M (\e^{2n-2} \, 
H_a + H_{h,k} ) \right| \leq \bar c_\kappa \, (r_\e^{2n+1} + 
\e^{4n-4}\, r_\e^{2-2n} ). 
\] 
\label{le:f-8.1} 
\end{lemma} 
{\bf Proof~:} The proof of these estimates can be obtained as in 
\cite{Are-Pac}. We briefly recall it here. First, we use the result 
of Proposition~\ref{pr:f-5.5bis} to estimate 
\begin{equation} 
\| H_{h,k} \|_{{\mathcal C}^{4, \alpha}_{3-2n} (M_{r_\e})} \leq 
c_\kappa \, r_\e^{2n+1} . \label{eq:jd} 
\end{equation} 
Now observe that $H_{h,k} =0$ in $M_{r_0}$ and Ê${\mathbb L}_0 \, 
H_{h,k} =0$ in each $B_{\ell,r_0/2} - B_{\ell ,r_\e}$, hence 
${\mathbb L}_M \, H_{h,k} = ({\mathbb L}_M -{\mathbb L}_0) \, 
H_{h,k}$ in this set. 

\medskip 

Now, we use the expansion (\ref{eq:normal}) which reflects the fact 
that the metric induced by $\omega_M$ is in each $B_{\ell , r_0}$ 
asymptotic to the Euclidean metric. This implies that there exists a 
constant $c >0$ such that 
\[ 
\| ({\mathbb L}_M - {\mathbb L}_0 ) \, \varphi \|_{{\mathcal C}^{0, 
\alpha}_{\delta - 4}(B_{\ell , r_0}^*)} \leq c \, \| \nabla2 \, 
\varphi \|_{{\mathcal C}^{2, \alpha}_{\delta - 4} (B_{\ell 
,r_0}^*)}. 
\] 
Collecting these, and taking into account the effect of the cutoff 
functions, we conclude that 
\[ 
\| {\mathbb L}_M \, H_{h,k} \|_{{\mathcal C}^{0, \alpha}_{\delta -4} 
(M_{r_\e})} \leq c_\kappa \, ( r_\e^{2n+1} + r_\e^{6-\delta}). 
\] 
The first estimate follows at once from the fact that $6 - \delta > 
2n+1$. 

\medskip 

The second estimate follows from the structure of the nonlinear 
operator $Q_M$ described in (\ref{eq:f-14}) together with the 
estimate 
\[ 
\| \nabla^{2} H_a \|_{{\mathcal C}^{2, \alpha}_{2-2n} (M_{r_\e})} 
\leq c . 
\] 
and also (\ref{eq:jd}). Roughly speaking this estimate reflects the 
fact the most important nonlinear terms in the expression of $Q_M$ 
are of the form $B_{\ell,j} (\nabla^{4} \, \varphi , \nabla^{2} 
\varphi)$ and $\tilde B_{\ell, j'} (\nabla^{3} \varphi , \nabla^{3} 
\, \varphi)$ provided the second derivatives of $\varphi$ remain 
bounded (which is precisely our case). The proof of the second 
estimate follows easily from these considerations. The third 
estimate is easy to obtain and left to the reader. \fdim 

\medskip 

We will also need the following result which will be useful for the 
contraction mapping argument we will apply ~: 
\begin{lemma} 
We fix $\delta Ê\in (4-2n , 5-2n)$ and $\tilde c_\kappa >0$. There 
exists a constant $\hat c_\kappa >0$ such that, for all $\varphi, 
\varphi' \in {\mathcal C}^{4, \alpha}_{\delta} (M_{r_\e})$ 
satisfying 
\[ 
\|\varphi \|_{{\mathcal C}^{4, \alpha}_{\delta }(M_{r_\e})\oplus 
{\mathcal D}} + \| \varphi' \|_{{\mathcal C}^{4, \alpha}_{\delta 
}(M_{r_\e}) \oplus {\mathcal D}} \leq \tilde c_\kappa \, ( 
r_\e^{2n+1} + \e^{4n-4} \, r_\e^{6 -4n -\delta }), 
\] 
we have 
\[ 
\begin{array}{rllll} 
\displaystyle Ê\| Q_M (\e^{2n-2} \, H_a + H_{h,k} + \varphi ) - Q_M 
(\e^{2n-2} \, H_a + ÊH_{h,k} + \varphi' ) \|_{{\mathcal C}^{0, 
\alpha}_{\delta -4}(M_{r_\e})} \qquad \qquad \\[3mm] 
\hfill \leq \hat c_\kappa \, \e^{2n-2} \, r_\e^{6-4n-\delta} \, \| 
\varphi - \varphi' \|_{{\mathcal C}^{4, \alpha}_{\delta 
}(M_{r_\e})\oplus {\mathcal D}} , 
\end{array} 
\] 
and we also have 
\[ 
\begin{array}{rllll} 
\displaystyle \int_{M(r_\e)} \left| Q_M (\e^{2n-2} \, H_a + H_{h,k} 
+ \varphi ) - 
Q_M (\e^{2n-2} \, H_a + ÊH_{h,k} + \varphi' ) \right| \qquad \qquad \\[3mm] 
\hfill \leq \hat c_\kappa \, \e^{2n-2} \, r_\e^{2-2n} \, \| \varphi 
- \varphi' \|_{{\mathcal C}^{4, \alpha}_{\delta }(M_{r_\e}) \oplus 
{\mathcal D} } . 
\end{array} 
\] \label{le:f-8.2} 
\end{lemma} 
The important observation is that 
\[ 
\lim_{\e \rightarrow 0} \e^{2n-2} \, r_\e^{6-4n-\delta} = \lim_{\e 
\rightarrow 0} \e^{2n-2} \, r_\e^{2-2n} = 0. 
\] 

Using the results of \S 5.2, we reduce the solvability of the 
equation (\ref{eq:f-22}) to a fixed point problem in 
\[ 
\{ \varphi Ê\in {\mathcal C}^{4, \alpha}_\delta (M_{r_\e}) \oplus 
{\mathcal D} \quad : \newline \quad Ê\|\varphi \|_{{\mathcal C}^{4, 
\alpha}_\delta (M_{r_\e}) \oplus {\mathcal D}} \leq c_\kappa \, 
(r_\e^{2n+1} + \e^{4n-4} \, r_\e^{6 -4n -\delta }) \} , 
\] 
where the constant $c_\kappa >0$ is fixed large enough. And in fact, 
using Lemma~\ref{le:f-8.1} and Lemma~\ref{le:f-8.2} it is easy to 
check that there exists a fixed point solution of (\ref{eq:f-22}) 
provided $\e$ is chosen small enough, say $\e \in (0, \e_0)$ by 
applying a fixed point theorem for contraction mapping. 

\medskip 

Hence we have obtained the~: 
\begin{prop} 
Assume that $n \geq 2$ and choose $\delta Ê\in (4-2n, 5-2n)$. Given 
$\kappa >0$, there exists $\breve{c}_\kappa >0$ and $\e_\kappa > 0$ 
such that, for all $\e \in (0, \e_\kappa)$, for all $h_\ell \in 
{\mathcal C}^{4, \alpha}(\del B_1)$ and all $k_\ell \in {\mathcal 
C}^{2, \alpha}(\del B_1)$ satisfying (\ref{eq:f-18}), there exists a 
\K\ form 
\[ 
\omega_{h,k} = \omega_M Ê+ i \, \del \, \bar \del \, \varphi_{h,k}, 
\] 
defined on $M_{r_\e}$, which has constant scalar curvature equal. 
Moreover, 
\[ 
\begin{array}{rllll} 
\| \varphi_{h,k} \, |_{B_{\ell ,2 r_\e} -B_{\ell ,r_\e}} (r_\e \, 
\cdot) - \e^{2n-2} \, a_\ell \, G_\ell (r_\e \, \cdot) - \e^{2n-2} 
\, b_\ell - H_{h_\ell, k_\ell}^o \|_{{\mathcal C}^{4, \alpha} ( 
B_2-B_1 )} \qquad \qquad \\[3mm] 
\hfill Ê\leq \breve{c}_\kappa \, (r_\e^{5} + \e^{4n-4} \, r_\e^{10 
-6n -\delta }), 
\end{array} 
\] 
and 
\[ 
|{\bf s} (\omega_{h,k}) - {\bf s} (\omega_M) + \e^{2n-2}\, a_0| \leq 
\breve{c}_\kappa \, (\e^{2n-1} + \e^{4n-4} \, r_\e^{2-2n}). 
\] 
\label{pr:f-8.1} 
\end{prop} 

Recall that $\varphi_\ell$ is the function which appears in 
(\ref{eq:normal}), Êthe expansion of the \K\ potential $\omega_M$ at 
$p_\ell$ so that 
\[ 
\omega_{h,k} = Êi \, \del Ê\, \bar \del \, ( \mbox{$\frac{1}{2}$} \, 
|z|^2 + \varphi_\ell + \varphi_{h,k}), 
\] 
in each $B_{\ell ,2 r_\e} -B_{\ell ,r_\e}$. We define, 
\[ 
\psi^o_\ell : = (\varphi_\ell + \varphi_{h,k}) (r_\e \, \cdot) + 
\e^{2n-2} \, b_\ell, 
\] 
when $n\geq 3$ and 
\[ 
\psi^o_\ell : = (\varphi_\ell + \varphi_{h,k}) (r_\e \, \cdot) + 
\e^{2n-2} \, b_\ell - \e^{2} \, a_\ell \, \log r_\e^{2} , 
\] 
when $n=2$. Observe that, locally, the \K\ potential is not uniquely 
defined and addition of constant does not modify the corresponding 
\K\ form. 

\medskip 

Thanks to the previous analysis, the following expansions are valid 
in each $B_2 -B_1$. 
\begin{lemma} 
The following expansions holds 
\begin{equation} \psi^o_\ell (z) = 
- a_\ell \, \e^{2n-2} \, r_\e^{4-2n} Ê\, |z|^{4-2n} + H^o_{h_\ell, 
k_\ell} (z) + {\mathcal O} (r_\e^{4}) + {\mathcal O}_\kappa 
(r_\e^{5} + \e^{4n-4} \, r_\e^{10 -6n -\delta }) , \label{eq:az-1} 
\end{equation} 
when $n \geq 3$ and 
\begin{equation} 
\psi^o_\ell (z) = a_\ell \, \e^{2} \, Ê\log Ê|z|^{2} + H^o_{h_\ell, 
k_\ell} (z) Ê+ {\mathcal O} (r_\e^{4}) + {\mathcal O}_\kappa ( 
r_\e^{5} + \e^{4} \, r_\e^{-2 -\delta }) , \label{eq:az-2} 
\end{equation} 
when $n =2$. \label{le:123} 
\end{lemma} 
Similar expansions are valid for the partial derivatives of these 
functions. The subscript $\kappa$ in ${\mathcal O}_\kappa (A)$ is 
meant to recall that, as $\e$ tends to $0$, this quantity is bounded 
by a constant, depending on $\kappa$, times $A$, while ${\mathcal 
O}(A)$ is bounded by a constant (independent of $\kappa$) times $A$. 

\medskip 

The key point is that 
\[ 
r_\e^{5} + \e^{4n-4} \, r_\e^{10 -6n -\delta } \ll r_\e4 
\] 
as $\e$ tends to $0$ (here we have to assume that $\delta \in (0, 
2/3)$, when $n=2$) and this implies that, in the expansion of the 
boundary data of the potential $\varphi_{h,k}$, the part 
corresponding to the solution of the nonlinear problem 
(\ref{eq:f-21}) is much smaller than $r_\e4$, the characteristic 
size of the boundary data $h$ and $k$. 

\medskip 

It should be clear that the \K\ form provided by 
Proposition~\ref{pr:f-8.1} depends continuously on the data (such as 
$h$ and $k$). 

\subsection{Perturbation of $\eta$} 

Given $\tilde a >0$, we now consider on $N$, the \K\ form 
\[ 
\eta_{\tilde a} : = (2^{2-n} \, \tilde a)^{\frac{1}{n-1}} \, \eta 
\] 
which, after a change of variable, can be expanded in $C_1$ as 
\begin{equation} 
\eta_{\tilde a } = i \, \del \, \bar \del \,\left( 
\mbox{$\frac{1}{2}$} \, |u|^2 - \tilde a \, |u|^{4-2n} + \tilde 
\varphi_{\tilde a} \right) , \label{eq:zdzd} \end{equation} when $n 
\geq 3$ and 
\begin{equation} 
\eta_{\tilde a} = i \, \del \, \bar \del \,\left( 
\mbox{$\frac{1}{2}$} \, |u|^2 + \tilde a \, \log |u|^2 + \tilde 
\varphi_{\tilde a} \right) , \label{eq:zdzdzd} \end{equation} when 
$n=2$, where the function $\varphi_{\tilde a} \in {\mathcal C}^{3, 
\alpha}_{3-2n} (C_1)$ enjoys properties similar to the one enjoyed 
by $\tilde \varphi$. Everything being uniform in $\tilde a $ as long 
as this parameter remains both bounded from above and bounded away 
from $0$. In fact we will need to choose $\tilde a$ so that 
\[ 
\frac{1}{2} \, \inf_\ell a_\ell < \tilde a < 2 \sup_{\ell} a_\ell , 
\] 
where the coefficients $a_\ell$ are the one which have been defined 
at the beginning of the previous section. 

\medskip 

As in the previous section, we consider the perturbed \K\ form 
\[ 
\tilde \eta = \eta_{\tilde a} + i \, \del \, \bar \del \, \varphi , 
\] 
The scalar curvature of $\tilde \eta$ is given by 
\[ 
{\bf s} (\tilde \eta) = Ê{\mathbb L}_N \, \varphi + Q_N (\varphi) , 
\] 
since the scalar curvature of $\eta$ is identically equal to $0$. 
Even though this is not explicit in the notations, the operators 
${\mathbb L}_N$ and $Q_N$ both depend on $t$. Again, the structure 
of the nonlinear operator $Q_{N}$ is quite complicated but, in 
$C_1$, it enjoys a decomposition similar to the ones described in 
(\ref{eq:f-14}) and the estimates (\ref{eq:f-15}) and 
(\ref{eq:f-16}) hold uniformly in each $C_1$ and in $\tilde a $, as 
long as $\tilde a $ remains bounded from above and bounded away from 
$0$. 

\medskip 

We define \begin{equation} R_\e : = Ê\frac{r_\e}{\e} \label{eq:RRe} 
\end{equation} 
Given $\tilde h \in {\mathcal C}^{4, \alpha}(\del B_1)$ and $\tilde 
k \in {\mathcal C}^{2, \alpha}(\del B_1)$ satisfying 
\begin{equation} 
\| \tilde h \|_{{\mathcal C}^{4, \alpha} (\del B_1)} + \| \tilde k 
\|_{{\mathcal C}^{2, \alpha} (\del B_1)} \leq \kappa \, R_\e^{3-2n}, 
\label{eq:f-23} 
\end{equation} 
where $\kappa >0$ will be fixed later on, we define 
\begin{equation} 
\tilde H_{\tilde h , \tilde k } : = Ê\tilde \chi \, (H^i_{\tilde h , 
\tilde k } (\cdot / R_\e) Ê- H^i_{\tilde h , \tilde k } (0) ) + 
H^i_{\tilde h , \tilde k} (0) , \label{eq:f-24} 
\end{equation} 
where $\tilde \chi$ are cutoff functions which are identically equal 
to $1$ in $C_2$ and identically equal to $0$ in $N_1$. 

\medskip 

We fix $\delta \in (0, 1)$. The space ${\mathcal C}^{4, 
\alpha}_{\delta} (N_{R_\e})$ (resp. $ {\mathcal C}^{0, 
\alpha}_{\delta-4} (N_{R_\e})$) is defined as the space of 
restrictions of functions in ${\mathcal C}^{4, \alpha}_{\delta} (N)$ 
(resp. ${\mathcal C}^{0, \alpha}_{\delta -4} (N)$) to $N_{R_\e}$, 
endowed with the induced norm. 

\medskip 

We would like to solve the equation 
\[ 
{\mathbb L}_{N} \, ( \tilde H_{\tilde h_\ell , \tilde k_\ell} + 
\varphi ) + Q_{N} ( \tilde H_{\tilde h_\ell , \tilde k_\ell} + 
\varphi ) = \e^{2} \, s 
\] 
where $s \in {\mathbb R}$ satisfies $|s|\leq 1 + |{\bf s} 
(\omega_M)|$. 

\medskip 

We consider an extension (linear) operator 
\[ 
\tilde {\mathcal E}_\e : {\mathcal C}^{0, \alpha}_{\delta-4} 
(N_{R_\e}) \longrightarrow {\mathcal C}^{0, \alpha}_{\delta-4} (N ) 
, 
\] 
such that $\tilde {\mathcal E}_\e \, \psi = \psi$ in $N_{R_\e}$, 
$\tilde {\mathcal E}_\e \, \psi$ is supported in $N_{2 R_\e }$ and 
$\| \tilde {\mathcal E}_\e\|\leq 2$. 

\medskip 

The equation we would like to solve can be rewritten as 
\begin{equation} 
{\mathbb L}_N \, \varphi Ê = \tilde {\mathcal E}_\e \left( \left( 
\e^{2} \, s - Ê{\mathbb L}_{N} \, \tilde H_{\tilde h , \tilde k } - 
Q_{N} ( \tilde H_{\tilde h , \tilde k} Ê+ \varphi ) 
\right)|_{N_{R_\e}} \right) . \label{eq:f-25} 
\end{equation} 

We prove the~: 
\begin{lemma} 
Assume that $\delta \in (0,1)$. There exists $c, \bar c_\kappa >0$ 
such that, for all $\e \in (0, \e_0)$ 
\[ 
\| {\mathbb L}_{N} \, \tilde H_{\tilde h , \tilde k } \|_{{\mathcal 
C}^{0, \alpha}_{\delta -4}(N_{R_\e})} \leq \bar c_\kappa \, 
R_\e^{2-2n} , \qquad \qquad \| \e2 \, s \|_{{\mathcal C}^{0, 
\alpha}_{\delta -4}(N_{R_\e})} \leq Êc \, R_\e^{3-2n -\delta} , 
\] 
and 
\[ 
\| Q_{N} (\tilde H_{\tilde h , \tilde k} ) \|_{{\mathcal C}^{0, 
\alpha}_{\delta -4}(N_{R_\e})} \leq \bar c_\kappa \, R_\e^{4-4n}. 
\] 
\label{le:f-8.3} 
\end{lemma} 
{\bf Proof~:} Again, the proof of these estimates can be obtained as 
in \cite{Are-Pac}. First, we use the result of 
Proposition~\ref{pr:f-5.5} to estimate 
\begin{equation} 
\| \nabla2 \, \tilde H_{\tilde h, \tilde k} \|_{{\mathcal C}^{2, 
\alpha}_{0} (N_{R_\e} - N_2)} \leq c_\kappa \, R_\e^{1-2n} \qquad 
\mbox{and} \qquad \| \nabla2 \, \tilde H_{\tilde h, \tilde k} 
\|_{{\mathcal C}^{2, \alpha}_{0} (N_2 - N_1)} \leq c_\kappa \, 
R_\e^{2-2n}. \label{eq:jddd} 
\end{equation} 
Now observe that ${\mathbb L}_0 \, \tilde H_{\tilde h, \tilde k} =0$ 
in each $C_1$, hence ${\mathbb L}_N \, \tilde H_{\tilde h, \tilde k} 
= ({\mathbb L}_N -{\mathbb L}_0) \, \tilde H_{\tilde h, \tilde k}$ 
in this set. 

\medskip 

Now, we use the expansion (\ref{eq:f-111}) or the expansion 
(\ref{eq:f-1111}) which reflect the fact that the metric induced by 
$\omega_N$ is, in each $C_1$ asymptotic to the Euclidean metric. 
This implies that there exists a constant $c >0$ such that 
\[ 
\| ({\mathbb L}_N - {\mathbb L}_0 ) \, \varphi \|_{{\mathcal C}^{0, 
\alpha}_{\delta - 4}(N_{R_\e} - N_2)} \leq c \, \| \nabla2 \, 
\varphi \|_{{\mathcal C}^{2, \alpha}_{\delta + 2n - 4} (N_{R_\e} 
-N_2)} \leq c \, \| \nabla2 \, \varphi \|_{{\mathcal C}^{2, 
\alpha}_{0} (N_{R_\e} -N_2)}. 
\] 
Collecting these, and taking into account the effect of the cutoff 
functions, we conclude that 
\[ 
\| {\mathbb L}_M \, H_{h,k} \|_{{\mathcal C}^{0, \alpha}_{\delta -4} 
(N_{R_\e})} \leq c_\kappa \, R_\e^{2-2n} . 
\] 
Observe that the main contribution comes from the effect of the 
cutoff function $\tilde \chi$. 

\medskip 

The second estimate is easy to derive. The third estimate follows 
from the structure of the nonlinear operator $Q_N$, which is similar 
to the one described in (\ref{eq:f-14}), together with 
(\ref{eq:jddd}). Again, the most important nonlinear terms in the 
expression of $Q_N$ are of the form $B_{j} (\nabla^{4} \, \varphi , 
\nabla^{2} \varphi)$ and $\tilde B_{j'} (\nabla^{3} \varphi , 
\nabla^{3} \, \varphi)$ provided the second derivatives of $\varphi$ 
remain bounded (which is precisely our case). The proof of the 
second estimate follows easily from these considerations and, again, 
the main contribution comes from the effect of the cutoff function 
$\tilde \chi$. \fdim 

\medskip 

We will also need the~: 
\begin{lemma} 
We fix $\delta \in (0,1)$ and $\tilde c >0$. There exists $\hat c > 
0$ such that, for all $\varphi, \varphi' \in {\mathcal C}^{4, 
\alpha}_\delta (N_{R_\e})$, satisfying 
\[ 
\| \varphi\|_{{\mathcal C}^{4, \alpha}_{\delta} (N_{R_\e})}+ 
\|\varphi'\|_{{\mathcal C}^{4, \alpha}_{\delta} (N_{R_\e})} \leq 
\tilde c \, R_\e^{3-2n -\delta} , 
\] 
we have 
\[ 
\| Q_{N} (\tilde H_{\tilde h , \tilde k} + \varphi) - Q_{N} (\tilde 
H_{\tilde h, \tilde k} + \varphi') \|_{{\mathcal C}^{0, 
\alpha}_{\delta -4}(N_{R_\e})} \leq \hat c \, R_\e^{3-2n-\delta} \, 
\| \varphi -\varphi'\|_{{\mathcal C}^{4, \alpha}_{\delta} 
(N_{R_\e})}. 
\] 
\label{le:f-8.4} 
\end{lemma} 

Using the result of \S 5.3, we reduce the solvability of the 
equation (\ref{eq:f-25}) to a fixed point problem in 
\[ 
\{ \varphi Ê\in {\mathcal C}^{4, \alpha}_\delta (N_{R_\e } ) \quad : 
\quad Ê\| \varphi \|_{{\mathcal C}^{4, \alpha}_\delta (N_{R_\e})} 
\leq \tilde c \, R_\e^{3-2n-\delta} \} , 
\] 
where the constant $\tilde c >0$ is fixed large enough. In fact, 
using Lemma~\ref{le:f-8.3} and Lemma~\ref{le:f-8.4} it is easy to 
check that there exists a fixed point solution of (\ref{eq:f-25}) 
provided $\e$ is chosen small enough, say $\e \in (0, \e_\kappa)$, 
which is obtained by applying a fixed point theorem for contraction 
mapping. 

\medskip 

We have obtained the~: 
\begin{prop} Given $\kappa >0$, $\tilde a >0$ and $s \in {\mathbb R}$, 
there exists $\e_\kappa >0$ such that, for all $\e \in (0, 
\e_\kappa)$, for all $\tilde h \in {\mathcal C}^{4, \alpha}(\del 
B_1)$ and $\tilde k \in {\mathcal C}^{2, \alpha}(\del B_1)$, there 
exists a \K\ form 
\[ 
\eta_{\tilde h , \tilde k , s , \tilde a} = \eta_{\tilde a} + i \, 
\del \, \bar \del \, \tilde \varphi_{\tilde h , \tilde k , s , 
\tilde a} , 
\] 
defined on $N (R_\e)$, which has constant scalar curvature. Moreover 
\[ 
\| \tilde \varphi_{\tilde h, \tilde k , s , \tilde a} |_{C_{R_\e } - 
C_{R_\e/2}} (R_\e \, \cdot ) - H^i_{\tilde h , \tilde k } 
\|_{{\mathcal C}^{4, \alpha} (B_1 -B_{1/2})} \leq c \, R_\e^{3-2n} , 
\] 
for some constant $c >0$ independent of $\kappa$. In addition, ${\bf 
s} \, (\tilde \omega_{\tilde h, \tilde k , s}) = \e^{2} \, s$. 
\label{pr:f-8.3} 
\end{prop} 

We now apply this result when 
\[ 
s = {\bf s}(\omega_{h,k}), 
\] 
is the scalar curvature of the metric defined in the previous 
section. The corresponding \K\ form will be denoted by $\eta_{\tilde 
h, \tilde k , \tilde a}$ even though it also depends on $h$ and $k$. 
Recall that $\tilde \varphi_{\tilde a}$ is the function which 
appears in (\ref{eq:zdzd}) and (\ref{eq:zdzdzd}), the expansion of 
the \K\ potential $\eta_{\tilde a}$ at $\infty$. Hence, the \K\ form 
$\eta_{\tilde h, \tilde k , \tilde a}$ has constant scalar curvature 
equal to $\e^{2} \, {\bf s} (\omega_{h,k})$ and, near $\infty$, it 
can be expanded as 
\[ 
\eta_{\tilde h, \tilde k , \tilde a} Ê= Êi \, \del \bar \del \, ( 
\mbox{$\frac{1}{2}$} \, |u|^2 - \tilde a \, |u|^{4-2n} + \varphi_t + 
\tilde \varphi_{\tilde h, \tilde k, s}), 
\] 
when $n \geq 3$ and 
\[ 
\eta_{\tilde h, \tilde k , \tilde a} Ê= Êi \, \del \bar \del \, ( 
\mbox{$\frac{1}{2}$} \, |u|^2 + \tilde a \, \log |u|^2 + \varphi_t + 
\tilde \varphi_{\tilde h, \tilde k, s}), 
\] 
when $n=2$. 

\medskip 

As in the previous section, we define, 
\[ 
\psi^i : = - \tilde a \, R_\e^{4-2n} \, |z|^{4-2n} + (\varphi_t + 
\tilde \varphi_{\tilde h, \tilde k, s}) ( R_\e \, \cdot ) , 
\] 
when $n\geq 3$ and 
\[ 
\psi^o : = \tilde a \, \log |z|^2 + (\varphi_t + \tilde 
\varphi_{\tilde h, \tilde k, s} ) ( R_\e \, \cdot), 
\] 
when $n=2$. Observe that, locally, the \K\ potential is not uniquely 
defined and addition of constant does not modify the corresponding 
\K\ form. 

\medskip 

Thanks to the previous analysis, the following expansions are valid 
in each $B_1 -B_{1/2}$. 
\begin{lemma} 
The following expansion holds 
\begin{equation} 
\psi^i (z) : = - \tilde a \, R_\e^{4-2n} \, |z|^{4-2n} + H^i_{\tilde 
h, \tilde k } (z) + {\mathcal O} (R_\e^{3-2n}) , \label{eq:az-3} 
\end{equation} 
when $n \geq 3$ and 
\begin{equation} 
\psi^i (z) = \tilde a Ê\, Ê\log Ê|z|^{2} + H^i_{\tilde h, \tilde k} 
(z) + {\mathcal O} (R_\e^{-1}) , \label{eq:az-4} 
\end{equation} 
when $n =2$. \label{le:1234} 
\end{lemma} Similar expansions are valid for the 
partial derivatives of these functions. 

\medskip 

The key point is that the quantity ${\mathcal O} (R_\e^{3-2n})$ is 
bounded by a constant, independent of $\kappa$ and of $\tilde a $, 
times $R_\e^{3-2n}$, as $\e$ tends to $0$. Again, it should be clear 
that the \K\ form provided by Proposition~\ref{pr:f-8.3} depends 
continuously on the data (such as $h, \tilde h , k $ and $\tilde 
k$). 

\section{Gluing the pieces together} 

We are now in a position to describe the connected sum construction. 
For all $\e >0$ Êsmall enough, we define a complex manifold $M_\e$ 
be removing small balls centered at the points $p_\ell$, $\ell= 1, 
\ldots , m $ and replacing them by properly rescaled versions of the 
$N$. We define 
\begin{equation} 
r_\e = Ê\e^{\frac{2n-1}{2n+1}}, Ê\qquad \mbox{and} \qquad R_\e : = 
\frac{r_\e}{\e} . \label{eq:f-9} 
\end{equation} 

By construction 
\[ 
M_\e : = M \sqcup _{{p_{1}, \e}} N \sqcup_{{p_{2},\e}} \dots \sqcup 
_{{p_{m}, \e}} N , 
\] 
is obtained by performing a connected sum of $M_{r_\e}$ (which has 
$m$ boundaries) with the truncated ALE spaces $N_{R_\e}, \ldots, 
N_{R_\e}$. This connected sum is obtained by identifying $\del 
B_{\ell \, r_\e}$ with $\del C_{R_\e}$ using 
\[ 
(z_1 , \ldots, z_n ) Ê= \e \, (u_1 , \ldots, u_n) , 
\] 
where $(z_1, \ldots, z_n)$ are coordinates in $B_\ell \, (r_0)$ and 
$(u_1, \ldots, u_n)$ are coordinates in $C_\ell \, (R_0)$. 

\medskip 

Next, a \K\ form is defined on $M_\e$ by gluing together the \K\ 
form $w_{h,k}$ define in \S 4.4 and the rescalled \K\ forms $\e^{2} 
\, \eta_{\tilde h_1, \tilde k_1, \tilde a_1}, \ldots, \e^{2} \, 
\eta_{\tilde h_m, \tilde k_m, \tilde a_m}$ defined in \S 4.5 on each 
copy of $N_{R_\e}$. 

\medskip 

It remains to explain to choose the data $h : = (h_1, \ldots, h_m)$, 
$k : = (k_0, \ldots, k_m)$ satisfying (\ref{eq:f-18}) Êand 
(\ref{eq:orth}), $\tilde h : = (\tilde h_1 , \ldots, \tilde h_m)$, 
$\tilde k : = (\tilde k_1 , \ldots, \tilde k_m)$ satisfying 
(\ref{eq:f-23}), and $\tilde a : = (\tilde a_1, \ldots, \tilde a_m)$ 
in such a way that, for each $\ell =1, \ldots, m$, the function 
$\psi_{\ell}^o$ defined in $B_{2}-B_{1}$ as in \S 4.4 on the one 
hand and the function $\e^{2} \, \psi_\ell^i$ defined in $B_{1} - 
B_{1/2}$ as in \S 4.5 with data $\tilde h_\ell, \tilde k_\ell$ and 
$\tilde a_\ell$ on the other hand, have their derivatives up to 
order $4$ which coincide on $\del B_{1}$. 

\medskip 

Now, assume that the functions $\psi_\ell^o$ and $\e^{2} \, 
\psi_\ell^i$ are ${\mathcal C}^{4, \alpha}$ functions which satisfy 
\begin{equation} 
\psi_\ell^o = \e^{2} \, \psi_\ell^i , \qquad \del_r \psi_\ell^o = 
\e^{2} \, \del_r \psi_\ell^i , \qquad \Delta_0 \psi_\ell^o = Ê\e^{2} 
\, \Delta_0 \psi_\ell^i, \qquad \del_r \Delta_0 \psi_\ell^o = \e^{2} 
\, \del_r \Delta_0 \psi_\ell^i, \label{eq:f-26} 
\end{equation} 
on $\del B_{1}$ where $r =|u|$. Consider polar coordinates and 
decompose 
\[ 
\Delta_0 = Ê\del_r2 + \frac{1}{r^{2}} \, \Delta_{S^{2n-1}} + 
\frac{2n-1}{r} \, \del_r. 
\] 
Using the decomposition of $\Delta_0$ in polar coordinates, it is 
clear that these identities guaranty that $\psi_\ell^o$ and $\e^{2} 
\, \psi_\ell^i$ have their partial derivatives up to order $3$ which 
coincide on $\del B_{1}$. 

\medskip 

Assuming for the moment that we have already found the construction 
data in such a way that (\ref{eq:f-26}) are satisfied, then the \K\ 
form 
\[ 
i \, \del \, \bar \del \, ( \mbox{$\frac{1}{2}$} \, r_\e^{2} \, 
|u|^2 + \psi^o_\ell ), 
\] 
defined in $B_{2} - B_{1}$ and the \K\ form 
\[ 
i \, \e^{2} \, \del \, \bar \del \, (\mbox{ $\frac{1}{2}$} \, 
R_\e^{2} \, |u|^2 + \psi^i_\ell), 
\] 
defined in $B_{1} - B_{1/2}$ have constant scalar curvature which 
are equal. This implies that the potential function defined by $\psi 
: = \psi_\ell^o$ in $B_{2} -B_{1}$ and $\psi : = \e^{2} \, 
\psi_\ell^i$ in $B_{1}- B_{1/2}$ is a weak solution of the nonlinear 
elliptic partial differential equation 
\[ 
{\bf s} \,\left( i \, \del\, \bar \del (\mbox{$\frac{1}{2}$} \, 
r_\e^{2} \, |u|^2 + \psi ) \right) = cte . 
\] 
in $B_2-B_{1/2}$. In addition $\psi$ is ${\mathcal C}^{3, \alpha}$ 
and smooth away from $\del B_{1}$. It then follows from elliptic 
regularity theory together with a bootstrap argument that this 
function is in fact smooth. Hence, by gluing the \K\ metrics 
$\omega_{h,k}$ and $\e^{2} \, \eta_{\tilde h_\ell, \tilde k_\ell , 
\tilde a_\ell}$ on the different pieces constituting $M_\e$, we have 
produced a \K\ metric on $M_\e$ which has constant scalar curvature. 
This will end the proof of the main Theorem. 

\medskip 

It remains to explain how to find the functions $h = (h_1, \ldots, 
h_m)$, $k = (k_1, \ldots, k_m)$, $\tilde h = Ê(\tilde h_1, \ldots, 
\tilde h_m)$, $\tilde k = Ê(\tilde k_1, \ldots, \tilde k_m)$ and the 
parameters $\tilde a Ê= (\tilde a_1, \ldots, \tilde a_m)$. To this 
aim, let us assume that $n \geq 3$ since only notational changes are 
needed to handle the general case $n = 2$. It follows from the 
result of Lemma~\ref{le:123} and Lemma~\ref{le:1234} that the 
following expansions hold 
\[ 
\begin{array}{rllll} 
\psi^o_\ell (z) & = & - a_\ell \, \e^{2n-2} \, r_\e^{4-2n} Ê\, 
|z|^{4-2n} + H^o_{h_\ell, k_\ell} (z) + {\mathcal O} (r_\e^{4}) 
\\[3mm] 
\e^{2} \, \psi^i_\ell (z) & = & - \tilde a_\ell \, \e^{2n-2} \, 
r_\e^{4-2n} \, |z|^{4-2n} + \e^{2} \, H^i_{\tilde h, \tilde k } (z) 
+ {\mathcal O} (r_\e^{4}) , 
\end{array} 
\] 
We change parameters and define the functions $h'_\ell, k'_\ell, 
\tilde h'_\ell$ and $\tilde k'_\ell$ by 
\[ 
\begin{array}{rllll} 
h'_\ell & : = & (\tilde a_\ell - a_\ell) \, r_\e^{4-2n} \, \e^{2n-2} 
+ h_\ell \\[3mm] 
k'_\ell & : = & 2 (4-2n) (\tilde a_\ell - a_\ell) \, \e^{2n-2} \, 
r_\e^{4-2n} Ê+ k_\ell 
\end{array} 
\] 
and also 
\[ 
\tilde h_\ell' := \e^{2} \, \tilde h_\ell \qquad \qquad \tilde 
k_\ell' := \e^{2} \, \tilde k_\ell 
\] 
Recall that the functions $k_\ell$ are assumed to satisfy 
(\ref{eq:orth}) while the function $k'_\ell$ is not assumed to 
satisfy such a constraint. The role of the scalar $\tilde a_\ell - 
a_\ell$ is precisely to recover this lost degree of freedom in the 
assignment of the boundary data. With these new variables, the 
expansions for both $\psi^o_\ell$ and $ \e2 \, \psi^i_\ell$ can now 
be written as 
\[ 
\begin{array}{rllll} 
\psi^o_\ell (z) & = & - \tilde a_\ell \, r_\e^{4-2n} \, |z|^{4-2n} + 
H^o_{h_\ell', k_\ell'} (z) + {\mathcal O} (r_\e^{4}) \\[3mm] 
\e2 \, \psi^i_\ell (z) & = & - \tilde a_\ell \, r_\e^{4-2n} \, 
|z|^{4-2n} + \, H^i_{\tilde h', \tilde k' } (z) + {\mathcal O} 
(r_\e^{4}) . 
\end{array} 
\] 
The data functions $h' :=(h'_1, \ldots, h'_m)$, $\tilde h_1', \ldots 
, \tilde h'_m)$ are assumed to be bounded by a constant $\kappa$ 
times $r_\e^{4}$ in $ {\mathcal C}^{4, \alpha} (\del B_1)$ and the 
data functions $k' : = (k'_1, \ldots, k'_m)$ and $\tilde k' : = 
(\tilde k_1', \ldots, \tilde k'_m)$ are assumed to be bounded by a 
constant $\kappa$ times $r_\e^{4}$ in $ {\mathcal C}^{2, \alpha} 
(\del B_1)$. The terms ${\mathcal O} (r_\e^{4})$ are nonlinear terms 
of $h',\tilde h',k', \tilde k'$ which are bounded by a constant 
(independent of $\kappa$ times $r_\e^{4}$ provided $\e$ is chosen 
small enough. 

\medskip 

We recall the following result~: 
\begin{lemma} 
\cite{Are-Pac} The mapping 
\[ 
\begin{array}{rclclll} 
\mathcal P :& \mathcal C^{4,\alpha}(\del B_1) \times \mathcal 
C^{2,\alpha}(\del B_1) & \longrightarrow & \mathcal 
C^{3,\alpha}(\del B_1 ) \times \mathcal C^{1,\alpha}(\del B_1 ) 
\\[3mm] 
& (h,k) Ê&\longmapsto Ê Ê& (\partial_{r} \, (H^o_{h, k}- ÊH^i_{h, 
k}), \partial_{r} \, \Delta_0 \, (H^o_{h, k}- ÊH^i_{h, k})) , 
\end{array} 
\] 
is an isomorphism. \label{le:f-8.5} 
\end{lemma} 

Using this Lemma, it is straightforward to check that the 
solvability of (\ref{eq:f-26}) reduces to a fixed point problem 
which can be written as 
\[ 
(h', Ê\tilde h' , Êk', \tilde k' ) = S_\e ( h', Ê\tilde h' , k', 
\tilde k) , 
\] 
where the nonlinear operator $S_\e$ satisfies 
\[ 
\| S_\e (h', Ê\tilde h' , Êk', \tilde k' ) \|_{({\mathcal C}^{4, 
\alpha})^{2m} \times ({\mathcal C}^{2, \alpha} )^{2m}} \leq c_0 \, 
r_\e^{4} , 
\] 
for some constant $c_0 >0$ which does not depend on $\kappa$, 
provided $\e$ is small enough. We finally choose 
\[ 
\kappa = 2 \, c_0 , 
\] 
and $\e \in (0, \e_\kappa)$, where $\e_\kappa$ is fixed small 
enough. We have therefore proved that $S_\e$ is a map from 
\[ 
A_\e : = \left\{ (h', Ê\tilde h' , Êk', \tilde k' ) \in ({\mathcal 
C}^{4, \alpha})^{2m} \times ({\mathcal C}^{2, \alpha})^{2m} \quad : 
\quad \| (h', \tilde h' , Êk', \tilde k' ) \|_{({\mathcal C}^{4, 
\alpha})^{2m} \times ({\mathcal C}^{2, \alpha})^{2m} } \leq \kappa 
\, r_\e^{4} Ê\right\} , 
\] 
into itself. This mapping is clearly continuous and if in addition 
it were compact, the application of Leray-Schauder's fixed point 
theorem would directly guaranty the existence of a fixed point. 
However, the nonlinear equations we have solved being fully 
nonlinear, there is no gain of regularity and hence 

\medskip 

To overcome this last difficulty, we define a family of smoothing 
mappings ${\mathcal D}_{q}$ Ê\cite{Alin} such that 
\begin{equation} 
\begin{array}{rllll} 
|| {\mathcal D}_q \, u ||_{{\mathcal C}^{r, \alpha}} & \leq & c \, 
||u ||_{{\mathcal C}^{r', \alpha'}} & \quad \mbox{for} \quad r + 
\alpha \leq r'+ \alpha' \\[3mm] 
|| {\mathcal D}_q \, u ||_{{\mathcal C}^{r, \alpha}} & \leq & c \, 
q^{r'+ \alpha' - r - \alpha} \, ||u ||_{{\mathcal C}^{r', \alpha'}} 
& \quad \mbox{for} \quad r + \alpha \geq r'+ \alpha'\\[3mm] 
|| u - {\mathcal D}_q \, u ||_{{\mathcal C}^{r, \alpha}} & \leq & 
c\, q^{r'+ \alpha' - r - \alpha} \, ||u ||_{{\mathcal C}^{r', 
\alpha'}} & \quad \mbox{for} \quad r + \alpha \leq r'+ \alpha'. 
\end{array} 
\label{pr} 
\end{equation} 
and replace the nonlinearities $S_\e$ by $S_{\e,q} : = {\mathcal 
D}_q \circ S_\e$. This time the corresponding operator $S_{\e,q}$ is 
compact and for $\e$ sufficiently small, maps the ball of radius 
$r_\e^{4}$ to itself. Hence it has a fixed point $(h'_q, \tilde 
h'_q, k'_q, \tilde k'_q)$ when $\e$ is small enough. Finally, the 
fixed points are bounded uniformly in $q \in (0,1)$, so for any 
fixed $\alpha' < \alpha$ we may extract a sequence $q_j \rightarrow 
0$ such that $(h'_{q_j}, \tilde h'_{q_j}, k'_{q_j}, \tilde 
k'_{q_j})$ converges in $({\mathcal C}^{4,\alpha'})^{2m} \times ( 
{\mathcal C}^{2,\alpha'})^{2m}$ to a fixed point of $S_\e$. This 
completes our proof of the existence the \K\ metric on $M_\e$ which 
has constant scalar curvature. 

\medskip 

Observe that the use of smoothing operator is not strictly necessary 
and the fixed point for $S_\e$ could have been obtained through the 
application of a fixed point theorem for contraction mappings, but 
this would have required more work in \S 4.4 and \S 4.5 to prove 
that $S_\e$ is a contraction for $\e$ small enough. 

\section{Understanding the constraints} 

In this section, we give the proof of Lemma 1.1 and Lemma 1.2. 

\medskip 

As usual, let us denote by $\xi_0 \equiv 1 , \xi_1 , \ldots, \xi_d$ 
the set of independent functions which span the kernel of ${\mathbb 
L}_M$. We assume that, for $j=1, \ldots, d$, the functions $\xi_j$ 
are normalized to have mean $0$ and to be mutually 
$L^{2}$-orthogonal. We keep the notations of the introduction and 
define the matrix ${\mathfrak M}$ as in (\ref{eq:d-8}) and the 
integer values functions ${\mathfrak C}_1$ and ${\mathfrak C}_2$ as 
in (\ref{eq:d-9}). 

\begin{lemma} 
Assume that $m \geq d$. Then the set of points $(p_1, \ldots, p_{m}) 
\in M^m$ such that ${\mathfrak C}_1 = d$ is open and dense in 
$M^{m}$. \label{le:f-9.9} 
\end{lemma} 
{\bf Proof~:} Observe that it is sufficient to consider the case 
where $m=d$, since increasing the number of columns of the matrix 
can only increase its rank ! 

\medskip 

The proof is by induction on $d$. When $d=0$ there is nothing to 
prove. Now, assume that the result if true for $d-1$ functions. We 
write the function $(p_1, \ldots, p_{d}) \longrightarrow \mbox{det} 
\, {\mathfrak M} (p_1, \ldots, p_{d})$ as a linear combination of 
functions depending on $p_{d}$, the coefficients of which are 
functions which depend on $p_1, \ldots, p_{d-1}$. In other words, we 
expand the determinant $\mbox{det} \, {\mathfrak M}$ with respect to 
the last row. We obtain 
\[ 
\mbox{det} \, {\mathfrak M} (p_1, \ldots, p_{d}) = M_{1}(p_1, 
\ldots, p_{d-1}) \, \xi_1(p_{d}) + \ldots + M_{d}(p_1, \ldots, 
p_{d_1}) \, \xi_d (p_{d}) 
\] 
where $M_{i}$ is Ê(up to factor $\pm 1$) equal to the determinant of 
the matrix $(\xi_a (p_b))_{a,b}$ where $a = 1, \ldots, \check i , 
\ldots ,d$ and $b = 1, \ldots, d-1$. By assumption, the set of 
$(p_1, \ldots, p_{d-1})\in M^{d-1}$ for which $(M_{1}, \ldots, 
M_{d}) \neq 0$ is open and dense in $M^d$. We fix $(p_1, \ldots, 
p_{d-1}) \in M^d$ such that $(M_{1}, \ldots M_{d}) \neq 0$. The 
functions $\xi_1, \ldots, \xi_{d}$ being independent, the set of 
$p_{d} \in M$ for which $\mbox{det} \, {\mathfrak M} \neq 0$ is open 
and non empty. Furthermore, unique continuation theorem for the 
solutions of ${\mathbb L}_M \, \xi =0$ implies that the the set of 
$p_{d} \in M$ for which $\mbox{det} \, {\mathfrak M} \neq 0$ is 
dense. This completes the proof of the result. \fdim 

\medskip 

The second condition for our construction to work asks for the 
existence $p_1, \ldots, p_m \in M$, with $m \geq d+1$, such that 
there exists a solution $a_1, \ldots, a_m >0$ to the system 
\[ 
\forall i=1, \ldots, d, \qquad \qquad \sum_{j=1}^m a_j \, \xi_i 
(p_j) = 0, 
\] 
This amounts to ask for the existence of points $p_1, \ldots, 
p_{m}\in M$ such that the kernel of ${\mathfrak M} : = {\mathfrak M} 
(p_1, \ldots, p_m)$ contains a vector whose entries are all 
positive, i.e. an element of $K^m_+$, the positive cone in ${\mathbb 
R}^n$. Equivalently, we have to prove that it is possible to find 
points $p_1, \ldots, p_m \in M$ such that the image of $^t 
{\mathfrak M}$ is included in a hyperplane of ${\mathbb R}^m$ whose 
normal belongs to $K^m_+$. 

\medskip 

Observe that, according to Lemma~\ref{le:f-9.9}, for all $m \geq 1$, 
there exists an open and dense set $U_m \subset M^m$ such that for 
all $(p_1, \ldots, p_m) \in U_m$, the kernel of ${\mathfrak M}$ is 
$(m-d)$-dimensional and in fact varies continuously as the points 
change in $U_m$. Equivalently, for any choice of the points in 
$U_m$, the image of $^t {\mathfrak M}$ is a $d$-dimensional subspace 
of ${\mathbb R}^m$ which also varies continuously as the points 
change in $U_m$. This clearly shows that the set of points $(p_1, 
\ldots, p_m) \in M^m$ such that ${\mathfrak C}_1 =d$ and ${\mathfrak 
C}_2 \neq 0$ is an open (probably empty !) set. Also, it should now 
be clear that once we have found $p_1, \ldots, p_m$ satisfying both 
conditions, then the conditions remain fulfilled after any 
adjunction of points to this list. 

\medskip 

Now, the condition that $\mbox{Im} \, (^t {\mathfrak M})$ is 
included in a hyperplane of ${\mathbb R}^m$ whose normal belongs to 
$K^m_+$ is equivalent to the requirement that $\mbox{Im} \, (^t 
{\mathfrak M})$ does not contain any vector of $K^m_+$, i.e. that 
all nonzero elements of $\mbox{Im} \, (^t {\mathfrak M})$ have 
entries which change sign. 

This being understood, for all $\Lambda : = (\lambda_1, \ldots, 
\lambda_d) \in S^{d-1}$, we define on $M$ the function 
\[ 
f_\Lambda Ê: = \sum_{j=1}^d \lambda_j \, \xi_j . 
\] 
The previous discussion can be summarized as follows : We have to 
prove that, it is possible to find $(p_1, \ldots, p_m) \in U_m$ such 
that, for all $\Lambda \in S^{d-1}$, the vector $(f_\Lambda (p_1), 
\ldots, f_\Lambda (p_m))$ has entries which change sign. 

\medskip 

As stated in the introduction, except in special cases, we have not 
been able to prove that the minimal number of points for which the 
above condition is satisfied is $d+1$ even though we suspect that 
this is the case. Nevertheless, we have the general result~: 
\begin{lemma} 
There exists $m_0 \geq d+1$ and, for all $m \geq m_0$, there exists 
an nonempty open set $V_m \subset M^m$ such that ${\mathfrak C}_2 
\neq 0$, for all $(p_1, \ldots, p_m) \in V_m$. 
\end{lemma} 
{\bf Proof~:} We keep the above notations. ÊPick $\Lambda \in 
S^{d-1}$. Then the function $p \in M \longrightarrow f_\Lambda (p)$ 
has mean $0$ (since it is a linear combination of the functions 
$\xi_1, \ldots, \xi_d$ which are Êare assumed to have mean $0$). 
Therefore it is possible to find $p_\Lambda , \tilde p_\Lambda \in 
M$ such that 
\[ 
f_\Lambda (p_\Lambda) < 0 < f_\Lambda ( \tilde p_\Lambda) 
\] 
Now, by continuity, we also have 
\[ 
f_{\Lambda'} (p) < 0 < f_{\Lambda'} (\tilde p) 
\] 
for all $\Lambda'$ in some open neighborhood $O_\Lambda$ of 
$\Lambda$ in $S^{d-1}$, all $p$ in some open neighborhood 
$o_\Lambda$ of $p_\Lambda$ in $M$ and all $\tilde p$ in some open 
neighborhood $\tilde o_\Lambda$ of $\tilde p_\Lambda$ in $M$. The 
sets $O_\Lambda$ constitute an open cover of $S^{d-1}$, and by 
compactness one can find a finite sub-cover 
\[ 
S^{d-1} = \cup_{j=1}^J O_{\Lambda_j} 
\] 
Given any $ (p_1, \ldots, p_J , \tilde p_1 , \ldots, \tilde p_J) \in 
\Pi_{j=1}^J \, o_{\Lambda_j} \times \Pi_{j=1}^J \, \tilde 
o_{\Lambda_j}$ and given any $\Lambda \in S^{d-1}$, it belongs to 
some $O_{\Lambda_{j}}$ and hence the $j$-th and the $(J+j)$-th 
entries of the vector 
\[ 
(f(p_1), \ldots, f(p_J) , f(\tilde p_1), \ldots, f(\tilde p_J)) 
\] 
do not have the same sign. Therefore, we have found $m = 2J$ points 
satisfying the required conditions. \fdim 

\medskip 

One can then define $m_0 \geq d+1$ to be the least number of points 
for which ${\mathfrak C}_1 =d$ and ${\mathfrak C}_2 \neq 0$. 

\section{Geometric interpretations of the constraints} 

Let us go back to a well known interpretation of the kernel of the 
operator ${\mathbb L}_M$ in terms of holomorphic vector fields. ÊWe 
recall some important results (which go back to Matsushima and 
Lichnerowicz) and emphasize only the part relevant to the geometric 
interpretation of the analysis we have carried out in the previous 
sections. We refer the readers to the papers \cite{lbs} and 
\cite{ls} and the book of Futaki \cite{Fu} for a complete 
introduction and other applications. Given a real valued function 
$\xi$ satisfying 
\[ 
{\mathbb L}_M Ê\xi =0 
\] 
and having called $X_\xi = \partial^{\#} \, \xi$ the $(1,0)$ part of 
the gradient of $\xi$, one can use (\ref{eq:f-3}) to conclude that 
$\bar{\partial} \, X_\xi = 0$, i.e. $X_\xi$ is a holomorphic vector 
field which vanishes somewhere on the manifold (since $\xi$ 
certainly has critical points !). Conversely, every Killing vector 
field which vanishes somewhere is the imaginary part of a 
holomorphic vector field of the form $X_\xi$ for some $\xi \in 
\mbox{Ker} \, {\mathbb L}_M$ (see for example \cite{Fu}, 
\cite{lbs},\cite{ls}). If in addition $(M,\omega_M)$ has constant 
scalar curvature, then, modulo the space of parallel holomorphic 
vector fields, every holomorphic vector field vanishing somewhere on 
$M$ arises as $X_\xi$ for some $\xi \in Ê\mbox{Ker} \, {\mathbb 
L}_M$ \cite{Fu}, \cite{lbs}, \cite{ls}. Of course the parallel part 
of the Lie algebra of holomorphic vector fields does not come in 
since the requirement of vanishing somewhere makes it trivial. We 
will therefore disregard (and in fact we have already done it in the 
statements of our results) these vector fields. 

\medskip 

The correspondence just described should not shadow some important 
differences in dealing with holomorphic vector fields and with 
function annihilating the operator ${\mathbb L}_M$. This appears 
clearly when comparing our sufficient conditions with the existence 
of holomorphic vector fields on the blow up manifold. The following 
is standard, but we include it for reader's convenience. 
\begin{prop} 
Let $X$ be a holomorphic vector field on $M$ and $p$ any point in 
$M$. Then $X$ lifts to a holomorphic vector field $\tilde{X}$ on 
$\mbox{Bl}_{p} \, M$ if and only if $X(p)=0$. Moreover any 
holomorphic vector field on $\mbox{Bl}_{p} \, M$ projects to a 
holomorphic vector field on $M$ vanishing at $p$. 
\end{prop} 
{\bf Proof~:} If $X$ lifts to $\tilde{X}$ on $\mbox{Bl}_{p} - M$, 
then the (holomorphic) flow associated to $\tilde{X}$ preserves 
$\mbox{Bl}_{p} \, M \setminus E$, where $E$ is the exceptional 
divisor, and $E$ is fixed. Therefore $X(p)=0$. Moreover being the 
blow up map an isomorphism away from the exceptional divisor, the 
same observation proves the last part of the proposition, since 
extending the projected vector field to zero at $p$ gives a 
continuous vector field and hence a holomorphic one. 

\medskip 

Conversely, we can look directly on a local chart in ${\mathbb 
C}^n$, with $p$ identified with the origin, since away from the 
origin the problem is trivial. Then 
\[ 
\mbox{Bl}_{0} \, {\mathbb C}^n = \{((z_1, \dots, z_n),[l_1,\dots, 
l_n]) \in {\mathbb C}^n \times {\mathbb P}^{n-1} \quad : \quad Êz_i 
\, l_j = z_j \, l_i, i,j=1, \dots , n \}. 
\] 

We consider the chart given by $u_i$ Êfor which the defining 
equations become $ z_1 = u_1u_n, \dots, z_{n-1}= u_{n-1}u_n, z_n = 
u_n$. In this chart, 
\[ 
\partial_{z_i} = \frac{1}{u_n}\partial_{u_1} \,\,\, i=1, \dots , n-1 
\] 
\[ 
\partial_{z_n} = -\frac{u_1}{u_n}\partial_{u_1}- \dots 
-\frac{u_{n-1}}{u_n}\partial_{u_{n-1}} + \partial_{u_n}. 
\] 

Hence a holomorphic vector field on ${\mathbb C}^n$, which can be 
written as $X= X^{\alpha}\partial_{z_{\alpha}}$ lifts to 
$\mbox{Bl}_{0} - {\mathbb C}^n - \{\mbox{Exceptional divisor}\}$ as 
the vector field 
\[ 
Y= \frac{X^{1}}{u_n}\partial_{u_1}+ \dots + \frac{X^{n-1}}{u_n} 
\partial_{u_{n-1}}+ X^n \, \left(-\frac{u_1}{u_n}\partial_{u_1}- \dots 
-\frac{u_{n-1}}{u_n}\partial_{u_{n-1}} + \partial_{u_n} \right). 
\] 

In order to extend $Y$ to the exceptional divisor we have to verify 
that the functions $\frac{X^{\alpha} - u_{\alpha}X^n}{u_n}$ are 
bounded for any $\alpha$ near $u_{\alpha}=0$. This clearly forces 
$X^{\alpha}(0) = 0$ for any $\alpha$, and therefore the vanishing of 
the vector field $X$ at the blow up point. Once this is verified the 
extension is continuous and hence holomorphic. \fdim 

\medskip 

Whether a holomorphic vector field (vanishing somewhere) lifts to a 
blow up depends only on its value at the point one is blowing up. 
Therefore, if the space of holomorphic vector fields with zeros has 
complex dimension $d$, blowing up $d$ points in generic position 
(generic in the sense described by the above proposition) one gets a 
manifold without such vector fields. In terms of the potentials used 
throughout this paper this translates in looking for $m$ points such 
that 
\[ 
\mbox{Rank} \, (\partial^{\#}(\xi_i)(p_j)) = d. 
\] 
Our first condition (${\mathfrak C}_1 = d$) is sensitive to the 
zeros of the bounded functions in $\mbox{Ker} \, {\mathbb L}_M$ with 
mean zero instead of their critical values. 

\medskip 

The second condition is certainly encoding a more subtle phenomenon. 
In fact, this second condition should be related to some suitable 
stability property of the blown up manifold. The fact that some 
positivity condition must hold is present in all known examples in 
different veins, and has been deeply investigated in the case of 
complex surfaces with zero scalar curvature by LeBrun-Singer 
\cite{lbs}, Rollin-Singer \cite{rs}, and for Del Pezzo surfaces by 
Rollin-Singer \cite{rs2}. In our construction it appears naturally 
from the analytical approach we developed, but a geometric 
interpretation is certainly worth seeking. 

\section{Examples} 

\subsection{The case of ${\mathbb P}^n$} 

A convenient way to study our problems on the blow up at points of 
the projective spaces is to look at ${\mathbb P}^n$ as the quotient 
of the unit sphere in Ê${\mathbb C}^{n+1}$ with complex coordinates 
$(z_{1}, \dots , z_{n+1})$ via the standard $S^{1}$-action given by 
the restriction of complex scalar multiplication. 

\medskip 

It is well known that the automorphism group of ${\mathbb P}^n$ is 
given by the projectivization of $\mbox{GL}(n+1, \mathbb C)$, whose 
complex dimension is $d = (n+1)^{2}-1$. We therefore seek for $d$ 
real functions whose $(1,0)$-part of the gradient generate the Lie 
algebra of the automorphism group as explained in the previous 
subsection. This can be done in two equivalent ways: either by 
explicit computation on the automorphism group, or by relying on the 
equivalence described in the previous section between this and the 
study of the kernel of the operator 
\[ 
{\mathbb L}_M = - \frac{1}{2}\, \Delta^{2}_M - \mbox{Ric}_M \cdot 
\nabla_M^{2} 
\] 
which for ${\mathbb P}^n$ with its Fubini-Study metric induced by 
the Hopf fibration becomes 
\[ 
{\mathbb L}_M = - \frac{1}{2}\, \Delta^{2}_{{\mathbb P}^n }- 
2(n+1)\Delta_{{\mathbb P}^n } = - \frac{1}{2}\, \Delta_{{\mathbb 
P}^n } (\, \Delta_{{\mathbb P}^n } + 4(n+1)). 
\] 
Our problem reduces to seeking a basis of functions with mean zero 
of the eigenspace of the Laplacian $\Delta_{{\mathbb P}^n }$ 
associated to the eigenvalue $4 \, (n+1)$ (i.e. the eigenspace of 
the Laplacian $\Delta_{S^{2n+1}}$ which are associated to the 
eigenvalue $4 \, (n+1)$ and are invariant under the $S1$ action), 
and this is clearly given by the $n^{2}+2n$ functions 
\[ 
\xi_{ab}= z_{a}\bar{z}_{b}+ z_{b}\bar{z}_{a}, \qquad \qquad Ê\hat 
\xi_{ab}= i Ê\, (z_{a}\bar{z}_{b} - z_{b}\bar{z}_{a}) 
\] 
for $1 \leq a < b \leq n+1$ and 
\[ 
\tilde \xi_{a} = Ê|z_{a}|^2 - |z_{a+1}|^2, 
\] 
for $a=1, \ldots, n$. Recall that we should add $\xi_0 \equiv 1$ to 
this list. 

\medskip 

Obviously any explicit calculation will be rather troublesome. It is 
hence very convenient (and giving best results) to introduce 
symmetries acting on the projective space in order to reduce as much 
as possible the elements of the kernel of ${\mathbb L}_{{\mathbb 
P}^n}$ which are invariant under these symmetries. 

\medskip 

{\bf Example 1:} Let us consider the group $G$ acting on ${\mathbb 
P}^n$ generated by the transformations 
\[ 
(z_{1},\ldots ,z_{n+1}) \longrightarrow ( \pm z_{1}, Ê\ldots, \pm 
z_{n+1}). 
\] 
Of course, the action of any element of the group on ${\mathbb 
C}^{n+1}$ maps the unit sphere into itself. The space of elements of 
the kernel of ${\mathbb L}_{{\mathbb P}^n}$ which are invariant 
under the action of the elements of $G$ is generated by $\{ 1, 
\tilde \xi_{1}, \ldots, \tilde \xi_{n}\}$. 

\medskip 

Let us fix the following set of blow up points 
\[ 
p_{1}= (1, 0 , \ldots ,0), \,\, p_{2}= (0,1, \ldots, 0), \, \ldots, 
\, p_{n+1}= (0, \ldots, 0,1). 
\] 
We find the matrix 
\[ 
{\mathfrak M} = \left( 
\begin{array}{ccccccc} 
1 & - 1 Ê& 0 Ê& \ldots & 0 Ê\\[2mm] 
0 & 1 & - 1 Ê& \ddots Ê& \vdots \\[2mm] 
\vdots Ê& Ê\ddots & \ddots & Ê\ddots Ê & 0 \\[2mm] 
0 & \ldots & 0& 1 & - 1 
\end{array} \right) 
\] 
whose rank is clearly equal to $n$ and has a one dimensional kernel 
spanned by the vector $(1, \ldots, 1)$ which has positive entries ! 
This proves that, working equivariantly with respect to the action 
of the group $G$, the blow up of ${\mathbb P}^n$ at the above $n+1$ 
points carries a constant scalar curvature \K\ form. 

\medskip 

We have then proved 
\begin{corol} 
The blow up of ${\mathbb P}^n$ at $p_{1}= (1,0, \ldots , 0), \ldots 
, p_{n+1}= (0, \ldots ,0, 1)$ has a constant scalar curvature \K\ 
metric. 
\end{corol} 
It is worth remarking that this result is optimal in the number of 
points to be blown up, since for fewer points the manifold would 
have nonreductive automorphisms group, and hence no \K\ metrics of 
constant scalar curvature by the Mathushima-Lichnerovicz 
obstruction. Another interesting aspect of this example is that the 
manifold obtained still has non trivial (in fact $n$ dimensional) 
automorphism group. The point is that the surviving automorphisms 
are precisely those which are not $G$-invariant. 

\medskip 

{\bf Exemple 2~:} Now, we still work equivariantly with respect to 
the action of the group $G$ defined above but we Êfix the following 
set of blow up points 
\[ 
p_{1}= (1, 0 , Ê\ldots ,0), \,\, p_{2}= (0,1, 0, Ê\ldots, 0), \, 
\ldots, \, p_{n}= (0, 0 , \ldots, 1 ,0) 
\] 
and 
\[ 
p_{n+1}= (0, \ldots, - \alpha , \beta) \qquad \qquad p_{n+2}= (0, 0 
, \ldots, \alpha , \beta). 
\] 
where $\alpha^{2}+ \beta^{2}=1$. This time, we find the matrix 
\[ 
{\mathfrak M} = \left( 
\begin{array}{ccccccc} 
1 & - 1 Ê& 0 Ê& \ldots & 0 & 0 Ê& 0\\[3mm] 
0 & 1 & - 1 Ê& \ldots Ê& 0 & 0 & 0\\[3mm] 
\vdots Ê& Ê& \ddots & \ddots Ê & Ê & Ê & \vdots \\[3mm] 
0 & \ldots & & 1 & -1 & - \alpha^{2} & - \alpha^{2}\\[3mm] 
0 & \ldots & & 0 & 1 & Ê\alpha^{2}-\beta^{2} & \alpha^{2}- \beta^{2} 
\end{array} \right) 
\] 
whose rank is again equal to $n$ and has a two dimensional kernel 
containing the vector 
\[ 
(1, \ldots, 1, \frac{1}{4\beta^{2}} , 1 - \frac{1}{2\beta^{2}} , 1- 
\frac{1}{2\beta^{2}}) 
\] 
which has positive entries, provided $0< \beta < \frac{1}{\sqrt{2}}$ 
! This proves that, working equivariantly with respect to the action 
of the group $G$, the blow up of ${\mathbb P}^n$ at the above $n+2$ 
points carries a constant scalar curvature \K\ form. 

\medskip 

It is an easy general observation that the addition of points to a 
list of points satisfying our sufficient conditions preserves these 
conditions satisfied. In particular we can add one point to the 
above lists, so, for example, the $G$-orbit of any point $p$ (which 
does not initially belong to the list) and keep the two conditions 
fulfilled. Observe that, for generic choice of the point $p$ the $G$ 
orbit of $p$ has $2^n$ points, so this substantially increases the 
number of points one can blow. 

\medskip 

However, one can add to one of the above lists $k$ points of the 
form 
\[ 
(0, \ldots, 0, \alpha , 0 , \ldots, \pm \beta , 0, \ldots ,0) 
\] 
where $\alpha^{2}+ \beta^{2}=1$, $\alpha \neq 0$ $\beta \neq 0$, so 
that the list of points remains invariant under the action of $G$. 
This clearly increases the number of blow up points by $2k$. Using 
this idea and starting from the list of points given in example 1, 
one shows that the blow up of ${\mathbb P}^n$ at $n+1 + 2 k$ points 
carries a constant scalar curvature \K\ form and Êstarting from the 
second list of points, one shows that the blow up of ${\mathbb P}^n$ 
at $n+2 + 2 k$ points carries a constant scalar curvature \K\ form. 
Therefore, we have obtained~: 
\begin{corol} 
The blow up of ${\mathbb P}^n$ at $m \geq n+1$ (special) points has 
a constant scalar curvature \K\ metric. 
\end{corol} 

\medskip 

{\bf Example 3~:} The important observation is that, so far we have 
worked equivariantly and the $m$ blow up points cannot be chosen 
into some open set of $M^m$. We now give a upper bound for the 
number $m_0$ which corresponds to the least number of points (larger 
than $d$) for which the conditions ${\mathfrak C}_1 = d$ Êand 
${\mathfrak C}_2 \neq 0$ are fulfilled. Given $\alpha, \beta \in 
{\mathbb R}-\{0\}$ satisfying $\alpha^{2} + \beta^{2} =1$ and 
$\alpha^{2}-\beta^{2} \neq 0$, we consider the following set of 
points~: 
\[ 
\begin{array}{rllll} 
p_{ij} & = & (0, \ldots, 0, \alpha , 0, \ldots, 0, \beta, 0\ldots ,0)\\[3mm] 
\tilde p_{ij} & = & (0, \ldots, 0, \alpha , 0, \ldots, 0, - \beta, 
0\ldots ,0)\\[3mm] 
\hat p_{ij} & = & (0, \ldots, 0, \beta , 0, \ldots, 0, i\, \alpha , 
0\ldots ,0)\\[3mm] 
\breve p_{ij} & = & (0, \ldots, 0, \beta , 0, \ldots, 0, -i \, 
\alpha , 0\ldots ,0) 
\end{array} 
\] 
where $1 \leq i<j \leq n+1$ correspond to the indices of the nonzero 
entries. There are exactly $2n(n+1)$ such points and it is easy to 
check that both conditions are fulfilled. Indeed, the first 
condition ${\mathfrak C}_1 =d$ is easy to check and left to the 
reader. Concerning the second condition, we show that $(1, \ldots, 
1)$ is in the kernel of ${\mathfrak M}$. Observe that this condition 
can be translated into the fact that for each $\xi = \xi_{ab}, 
\hat\xi_{ab}$ or $\tilde \xi_{a}$ 
\[ 
\sum_p \xi (p) = 0 
\] 
where summation over $p \in \{ p_{ij}, \tilde p_{ij}, \hat p_{ij}, 
\breve p_{ij} \quad : \quad i,j\}$ is understood. Observe that it is 
enough to check that the formula holds for the functions $z_a \, 
\bar z_b$ and the functions $|z_{a+1}|^2- |z_a|^2$. It is now easy 
too check that ${\mathfrak C}_2 \neq 0$. 

\medskip 

Therefore, we have obtained the~: 
\begin{corol} 
When $M = {\mathbb P}^n$, we have the estimate $m_0\leq 2 n (n+1)$. 
\end{corol} 

{\bf Example 4~:} For general \K\ manifolds with space of 
holomorphic vector fields of dimension $d$, one needs to blow up $d$ 
points in general position to have a manifold without holomorphic 
vector fields. In this respect ${\mathbb P}^n$ is very special since 
it is easy to observe that $n+2$ points suffice provided they form a 
so called projective frame, namely any choice of $n+1$ of them are 
linearly independent in ${\mathbb C}^{n+1}$ (such sets of points are 
also often called {\em{in generic position with respect to 
hyperplanes}}). The freedom of choices of projective frames ranges 
clearly in an open and dense subset of $({\mathbb P}^n)^{n+2}$. 

\medskip 

Suppose that we find a projective frame $p_1, \dots , p_{n+2}$ for 
which we can prove that $Bl_{p_1, \dots p_{n+2}}{\mathbb P}^n$ has a 
\K\ constant scalar curvature metric necessarely, being $n+2 < d$, 
using some equivariant construction. ÊThe first simple observation 
is that $Bl_{p_1, \dots , p_{n+2}}{\mathbb P}^n$ is biholomorphic to 
$Bl_{q_1, \dots ,q_{n+2}}{\mathbb P}^n$ for any other projective 
frame $q_1, \dots , q_{n+2}$. The manifold obtained is also with at 
most discrete automorphisms and we can then apply to it the results 
of \cite{Are-Pac} to blow up any other set of points and still get 
\K\ constant scalar curvature metrics. 

\medskip 

For all these reasons we now seek for a projective frame for which 
some equivariant construction works. To this aim, consider the group 
$G$ of permutations on the $n+1$ affine coordinates and the points 
\[ 
p_{1}= \frac{1}{\sqrt{n+1}}(1,\dots,1), \qquad \qquad p_{2}= 
\frac{1}{\sqrt{n+\alpha^{2}}}(1, ,\dots,1,- \alpha) 
\] 
with $\alpha>0$. ÊThe equivariant kernel is spanned by $\{ 1, 
\sum_{a<b}\xi_{ab}\}$ and therefore the matrix we obtain is 
\[ 
{\mathfrak M} = \left( 
\begin{array}{ccccccc} 
n & \frac{n(n-1-2\alpha)}{n+\alpha^{2}} 
\end{array} \right) 
\] 
Our sufficient conditions are then fulfilled for $\alpha > \frac{n-1}{2}$. 
It is easy to see that the points given by the orbits of $p_1$ 
(which is fixed), and $p_2$ (which consists of $n+1$ elements) form 
a projective frame and hence the resulting manifold is without 
holomorphic vector fields. 

\begin{corol} 
Given $m \geq n+2$ points in ${\mathbb P}^n$ which contain a 
projective frame, then the blow up of ${\mathbb P}^n$ at these $m$ 
points has a constant scalar curvature \K\ metric. Moreover the 
$n+2$ of these points which form a projective frame can move in an 
open set of $M^{n+2}$ and the remaining $m-n-2$ are arbitrary. 
\end{corol} 

This example naturally rises an important comment. For $3\leq m\leq 
8$ Tian's solution of the Calabi Conjecture tells us that 
$\mbox{Bl}_{p_1, \ldots, p_m} \, Ê{\mathbb P}^{2}$ admits \K 
-Einstein metrics (of course in the classes Ê$[\omega_{FS}] - 
[E_{1}] - \cdots - [E_{m}]$) as long as no three Êcollinear points 
are blown up, no $5$ of them lie on a quadric and Ê$8$ on a cubic. 
These conditions of course trivially prevent the Êblown up manifold 
to have a \K -Einstein metric. Moreover it is well known Êthat, 
under these assumptions, for $4\leq m \leq 8$, $\mbox{Bl}_{p_1, 
\ldots, p_m} \, {\mathbb P}^{2}$ has only discrete automorphisms. 
Therefore, as noted in \cite{Are-Pac}, we can use these metrics for 
successive blow ups and apply the main result in \cite{Are-Pac} to 
represent the classes $[\omega_{FS}] - [E_{1}] - \cdots - [E_{m}] - 
\epsilon [E_{m+1}]$ with metrics of constant scalar curvature. What 
we have gained with the present construction is firstly to represent 
also the classes $[\omega_{FS}] - \epsilon [(E_{1} + \cdots + 
E_{m})]$ with canonical positively curved metrics, therefore getting 
another open set (by perturbation arguments) in the \K\ cone. 
Applying their technique on parabolically stable bundles over 
riemann surfaces Rollin-Singer proved the same result as ours for 
$m\geq 7$. 

\medskip 

Moreover for $m>4$ the open set of points for which our construction 
works is larger than the above mentioned conditions for the 
existence of a \K\ -Einstein metric. For example once four points 
forming a projective frame are blown up, we can add a fifth aligned 
to two of the previous ones still getting a canonical metric. 

\medskip 

{\bf Example 5~:} The case of $M = {\mathbb P}^1 \times M_0$. The 
type of manifolds we now treat are of particular interest since they 
have been a central object of study in the papers of LeBrun-Singer, 
Kim-LeBrun-Pontecorvo up to the recent works of Rollin-Singer, when 
$M_0$ is taken to be a Riemann surface. The point is that a complete 
understanding of these examples leads {\it via} algebraic geometric 
techniques to the relation with stability of rank two vector bundles 
over Riemann surfaces. By understanding via Êa different approach 
these models in our more general setting, we hope to give a tool to 
the study of similar approach to higher rank vector bundles over any 
\K\ constant scalar curvature manifold. 

\medskip 

We assume throughout this example that $M_0$ is a \K\ manifold of 
any dimension with no nonvanishing holomorphic vector fields. 
Moreover on ${\mathbb P}^{1} \times M_0$ we consider the product 
metric, where ${\mathbb P}^{1}$ is endowed with the Fubini-Study 
metric normalized as in the previous subsection though for any 
positive value of the scalar curvature the results we obtain would 
hold verbatim. 

\medskip 

With these conventions the bounded kernel of the operator 
\[ 
{\mathbb L}_M = - \frac{1}{2}\, \Delta^{2}_M - \mbox{Ric}_M \cdot 
\nabla_M^{2} 
\] 
is naturally identified with $\xi_{12}$, $\hat\xi_{12}$ and $\tilde 
\xi_{a}$ the functions on ${\mathbb P}^{1}$ which have been 
described above. 

\medskip 

We then look at the group generated by 
\[ 
(z_{1},z_{2}) \rightarrow (z_{2}, -z_{1}), 
\] 
which reduces the invariant kernel to be generated by $\{1, 
\hat\xi_{12}\}$. Let us choose the points $p_{1}= 
\frac{1}{\sqrt{2}}(i,1)$ and $p_{1}= 
\frac{1}{\sqrt{2}}(1,i)$. The first condition then requires the 
matrix 
\[ 
{\mathfrak M} = \left( 
\begin{array}{ccccccc} 
-1 & 1 
\end{array} \right) 
\] 
Obviously the two conditions are fulfilled. Since the $G$-orbits of 
the points chosen are projectively the points themselves, we have 
proved~: 
\begin{corol} 
Given any two points $q_{1}$ and $q_{2}$ in $M_0$ (possibly 
coinciding), the blow up of ${\mathbb P}^1 \times M_0$ at $(p_{1}, 
q_{1})$, $(p_{2}, q_{2})$ has a constant scalar curvature \K\ 
metric. 
\end{corol} 
A stronger version of the above Corollary was proved by 
LeBrun-Singer when $M_0$ is a Riemann surface of genus at least $1$ 
in \cite{lbs} (in the case $M_0$ is a torus one must recall that all 
holomorphic vector fields are parallel so they do not interfere in 
the analysis of the operator ${\mathbb L}_M$). In fact they proved 
that the metric on the blow up can be chosen to have zero scalar 
curvature. For $M_0 = \Sigma_g$ we have gained the freedom of 
assigning any sign to the constant scalar curvature abtained. 

\medskip 

{\bf Example 6~:} The case of ${\mathbb P}^{n} \times {\mathbb 
P}^{m}$. The results for this case easily follow from the 
corresponding results for the previous examples. 

\medskip 

The case of ${\mathbb P}^{1} \times {\mathbb P}^{1}$ falls directly 
in the previous discussion since the blow up of this manifold at $m 
\geq 1$ points is biholomorphic to the blow up at $m+1$ suitably 
choosen points of ${\mathbb P}^{2}$. By the previous calculation we 
already know that for $m=2$ our problem has a positive solution, and 
that this is the least number of points for which this can happen. 
Moreover we have already seen that for $3$ points onwards (suitably 
chosen in an open set) the manifold has only discrete automorphisms, 
therefore iteration can start. 

\medskip 

Concerning the higher dimensional examples, for sake of simplicity 
we illustrate the case of ${\mathbb P}^{1} \times {\mathbb P}^{2}$, 
but everything can be easily extended to any of these products. ÊLet 
us consider affine coordinates $(z_{1}, z_2)$ on ${\mathbb P}^{1}$ 
and $(w_1,w_2,w_3)$ on ${\mathbb P}^{2}$ and let us denote by 
$\xi_{ab}, \hat \xi_{ab}, \tilde \xi_a$ the functions elements of 
the kernel on ${\mathbb P}^{1}$ as described above and by $\chi_{ab} 
, \hat \chi_{ab}, \tilde \chi_a$ the functions elements of the 
kernel on ${\mathbb P}^{2}$. To get the optimal result we have to 
work equivariantly and we consider the diagonal action used in the 
first example. ÊThe equivariant bounded kernel is then spanned by 
the union of the generators of each single factor, hence by $\tilde 
\xi_{1}$, $\tilde \chi_{1}$ and $\tilde \chi_{2}$. 

\medskip 

Let us choose the points 
$$ p_1=((1,0),(0,1,0)) \qquad Êp_2=((0,1),(1,0,0))$$ 
$$p_3=((1,0),(0,0,1))\qquad Êp_4 = 
(\frac{1}{\sqrt{5}}(1,2),(\frac{1}{\sqrt{3}}(1,1,1))\,\, .$$ 
The matrix Êis then 
\[ 
{\mathfrak M} = \left( 
\begin{array}{rrrrrr} 
-1 & 1 & 1 & -\frac{3}{5}\\ 
-1 & 1& 0 & 0& \\ 
1 & 0 & -1 & 0 
\end{array} \right) 
\] 
which is clearly of rank $3$ and whose kernel contains elements with 
positive entries. Therefore the blow up has a \K\ constant scalar 
curvature metric. It is moreover not difficult to see that the blow 
up of ${\mathbb P}^{1} \times {\mathbb P}^{2}$ has no holomorphic 
vector fields and therefore can be used both for moving these points 
in some open set and for starting iterations of further blow ups. 
Note that this examples generalizes immediately to any dimension in 
the following 

\begin{corol} Suppose $n\leq m$. 
The blow up of ${\mathbb P}^{n} \times {\mathbb P}^{m}$ at $k \geq 
m+2$ special points has a constant scalar curvature \K\ metric. 
Moreover, if the projections of these points to the two factors 
contain projective frames, these points can move in an open set of 
$({\mathbb P}^{n} \times {\mathbb P}^{m})^{m+2}$ and the other ones 
are arbitrary. 
\end{corol}

\end{document}